\documentclass[12pt]{article}
\usepackage{amsmath,amssymb,amsthm, amscd}

\begin{document}
\begin{center} 
{\bf Equivalence of symplectic singularities} 
\end{center} 
\vspace{0.4cm}

\begin{center}
Yoshinori Namikawa
\end{center} 
\vspace{0.4cm}

\begin{abstract}  
After introducing an equivalence problem for symplectic singularities, 
we formulate an algebraic version of such a problem. 
Let $X$ be an affine normal variety with a $\mathbf{C}^*$-action 
having only positive weights. Assume that the regular part $X_{reg}$ 
of $X$ admits an algebraic symplectic 2-form $\omega$ with weight $l$. 
Our main theorem asserts that any algebraic symplectic 2-form $\omega'$ on $X_{reg}$ 
of weight $l$ is equivalent to $\omega$ up to $\mathbf{C}^*$-equivariant 
automorphism of $X$ if $l \ne 0$.       
When $l = 0$ we have a counter-example to this statement. 
In the latter half of the article, we discuss the equivalence 
problem {\em up to constant}. We associate to $X$ a projective variety 
$\mathbf{P}(X)$ and prove that $\mathbf{P}(X)$ has a contact orbifold 
structure. Moreover, when $X$ has canonical singularities, the contact 
orbifold structure is rigid under a small deformation.  
The equivalence problem is then reduced to the uniqueness of the contact 
structures.  
In most examples the symplectic structures turn out to be unique up to 
constant with very few exceptions. In the final section we pose a splitting 
conjecture for symplectic singularities.   
\end{abstract}

\begin{center}
{\bf Introduction}
\end{center}
 
Assume that $X$ is a germ of a normal complex space whose regular locus $X_{reg}$ 
admits a holomorphic symplectic 2-form $\omega$.  
Two such pairs $(X,\omega)$ and $(X', \omega')$ 
are {\em equivalent} if there is an isomorphism $\phi : X \to X'$ such that 
$\omega = \phi^*(\omega')$.  They are not, a priori, equivalent 
even if their underlying complex analytic structure are equivalent. The Darboux theorem  
asserts that any holomorphic symplectic structure on $(\mathbf{C}^{2n},0)$ is 
equivalent to the standard one $dx_1 \wedge dx_2 + ... + dx_{2n-1} \wedge dx_{2n}$.
A general theme of this article is such an equivalence problem for a singular 
space. 

The Darboux theorem is naturally extended to a symplectic quotient singularity (cf. 
Proposition (1.1)). An essential idea for proving the Darboux theorem is due to Moser [Mo] and 
it seems rather difficult to develop this method for an arbitrary singular space.  

In this article we formulate algebraic versions of the equivalence 
problem. Namely we start with a normal affine variety $X$ of 
dimension $2d$ with a $\mathbf{C}^*$-action. Assume that $0 \in X$ is a unique fixed point of the 
$\mathbf{C}^*$-action with positive weights. More precisely, the cotangent space $m_{X,0}/m^2_{X,0}$ of $0 \in X$ has 
only positive weights with respect to the $\mathbf{C}^*$-action or equivalently, the affine ring $R$ of $X$ 
is positively graded: $\oplus_{i \ge 0}R_i$ with $R_0 = \mathbf{C}$. We call such a 
$\mathbf{C}^*$-action {\em good}.   
Let $\omega$ be an algebraic symplectic 
2-form on $X_{reg}$ with weight $l$. If we represent the $\mathbf{C}^*$-action by 
the family $\{\phi_t\}_{t \in \mathbf{C}^*}$ of automorphisms of $X$, then $\phi_t^*(\omega) = t^l \cdot \omega$. 
If we change the $\mathbf{C}^*$-action of $X$, $l$ may possibly change. But the positivity of $l$ 
reflects the property of $X$ itself. In fact, if $X$ has canonical singularities, then 
$l$ must be positive (Lemma (2.2)). Conversely, when $l >0$, one can show that $X$ has canonical singularities  
under the assumption that $X$ has an isolated singularity (Lemma (2.4)).   
 
Let us consider two such pairs $(X, \omega)$ and $(X', \omega')$ with the same weight 
$l$. They are called 
equivalent if there is an $\mathbf{C}^*$-equivariant isomorphism $\phi: X \cong 
X'$ such that $\omega = \phi^*(\omega')$. In particular, if $X = X'$, then 
$\omega$ and $\omega'$ are called equivalent symplectic structures on $X$. 
Our main result is: 
\vspace{0.2cm}

{\bf Theorem (3.1)} {\em Let $(X, \omega)$ be a pair of a normal affine variety $X$  
with a good $\mathbf{C}^*$-action and an algebraic symplectic 2-form $\omega$ 
on $X_{reg}$ with weight $l \ne 0$. Then  $\omega$ is a unique  
symplectic structure with weight $l$ up to equivalence. }      
\vspace{0.2cm}

If we drop the assumption $l \ne 0$, then the result does not hold. 
We have a counter-example when $l = 0$ (Remark (3.3)).  

Affine symplectic varieties are constructed in several different manners such as 
nilpotent orbit closures of a complex simple Lie algebra (cf. [C-M]), Slodowy 
slices to such orbits (cf. [Sl]) and the symplectic (or hyperK\"{a}hler) reductions. 
Note that these examples naturally come up with $\mathbf{C}^*$-actions.
It often happens that the same $\mathbf{C}^*$-variety appears in different constructions.  
But Theorem (3.1) asserts that the symplectic structures on the same $\mathbf{C}^*$-variety 
are unique if they have the same weight. We shall explain how Theorem (3.1) is applied to explicit examples.    
 
Let $\mathbf{g}$ be a complex simple Lie algebra and let $S \subset \mathbf{g}$ 
be a Slodowy slice to a nilpotent orbit $O$ of $\mathbf{g}$. Let 
$\mathbf{h}$ be a Cartan subalgebra and $W$ the associated Weyl group of $\mathbf{g}$.      
We denote by $\chi: S \to \mathbf{h}/W$ the adjoint quotient map restricted to 
$S$. We write $S_0$ for the central fibre $\chi^{-1}(0)$ of $\chi$. 
It is known that $S_0$ admits a (Kostant-Kirillov) symplectic structure $\omega$ together with a $\mathbf{C}^*$-action 
such that $\omega$ has weight 2. 
\vspace{0.15cm}

{\bf Example 1}. (i) Let $\mathbf{g}$ be the simple Lie algebra of type $B_n$ (resp. 
$C_n$, $F_4$ or $G_2$). Let $\mathbf{g}'$ be the simple Lie algebra of type 
$A_{2n-1}$ (resp. $D_{n+1}$, $E_6$ or $D_4$). 
Consider the Slodowy slices $S$ and $S'$ respectively for the subregular orbits 
of $\mathbf{g}$ and $\mathbf{g}'$. Then both $S_0$ and $S'_0$ have a Du Val singularity  
of type $A_{2n-1}$ (resp. $D_{n+1}$, $E_6$ or $D_4$). Moreover they are isomorphic as 
$\mathbf{C}^*$-varieties ([Sl], 7.4, Proposition 2, and 8.3, Proposition 2). 
According to Theorem (3.1)\footnote{In this case one can check 
easily that they are symplectic equivalent without Theorem (3.1).} we see that $(S_0, \omega)$ 
and $(S'_0, \omega')$ are also equivalent as symplectic varieties. 
This fact has a nice application to the study of Poisson deformations of 
$(S_0, \omega)$.  
As is observed in [LNS], the Poisson deformation $S \to \mathbf{h}/W$ 
of $(S_0, \omega)$ is not the universal one. But, since $(S_0, \omega) \cong (S'_0, \omega')$, 
one can regard $S' \to \mathbf{h}'/W'$ as a Poisson deformation of $(S_0, \omega)$. 
Since $\mathbf{g}'$ is simply-laced, this turns out to be universal. 

(ii)  Let $\mathbf{g}$ be the simple Lie algebra of type $G_2$ and 
let $S$ be a Slodowy slice to the 8-dimensional nilpotent orbit of $\mathbf{g}$. 
Let $\mathbf{g}'$ be the simple Lie algebra of type $C_3$  and let $S'$ be a 
Slodowy slice to the nilpotent orbit of $\mathbf{g}'$ of Jordan type $[4,1^2]$. 
Then $S_0$ and $S'_0$ are isomorphic to the 4-dimensional quasihomogeneous hypersurface 
$$ X:= \{(a,b,x,y,z) \in \mathbf{C}^5; a^2x + 2aby + b^2z + (xz-y^2)^3 = 0\} $$
as $\mathbf{C}^*$-varieties ([LNS], Section 10).  
Then $(S_0, \omega)$ and $(S'_0, \omega')$ are equivalent as symplectic varieties by 
Theorem (3.1). As in (i), $S \to \mathbf{h}/W$ does not give the universal Poisson 
deformation of $(S_0, \omega)$. But $S' \to \mathbf{h}'/W'$ is the universal Poisson 
deformation of $(S_0, \omega)$. 
\vspace{0.15cm}

{\bf Example 2}([L-N-S-vS]): {\em Quasihomogeneous symplectic hypersurfaces} 

At this moment we know two kinds of quasihomogeneous symplectic hypersurfaces.  
The first one is a series of examples $X_n$, $(n \geq 2)$ of dimension 4:   
$$ X_n:= \{(a,b,x,y,z) \in \mathbf{C}^5; a^2x + 2aby + b^2z + (xz-y^2)^n = 0\}. $$ 

The second one is a 6-dimensional example. For details on this example, see 
[LNS], Section 10. 

One can put (homogeneous) symplectic structures on them in several different ways. 
 
(a) Originally these examples were found as the central fibres $S_0$ of the Slodowy slices $S$ 
to certain nilpotent orbits of 
$\mathbf{g}$. The $X_n$ is the $S_0$ for the nilpotent orbit $O_{[2n-2, 1^2]}$ of 
$sp(2n)$ and the 6-dimensional example is the $S_0$ for the (unique) 6-dimensional 
nilpotent orbit of $G_2$. A Slodowy slice has a $\mathbf{C}^*$-action 
and admits a symplectic structure of weight $2$.   

(b) Let $V$ be an even dimensional representation of $sl_2$. One can put  
a Poisson structure on $A:= \mathbf{C}[sl_2 \oplus V]$ by using the Lie bracket 
of $sl_2$, the $sl_2$-representation $V$ and an $sl_2$-equivariant map 
$\varphi: \wedge^2 V \to \mathbf{C}[sl_2]$.  More precisely, for  
$x + v$, $y + w \in sl_2 \oplus V$, we define 
$$ \{x + v, y + w\}:= [x,y] + \varphi (v,w) + (x\cdot w + y\cdot v)$$ 
and extend this bracket to a Poisson structure on $A$ by the Leibniz rule.    

Take as $V$ the standard 2-dimensional representation 
and take as $\varphi$ the $(n-1)$-th power $\Delta^{n-1}$ of the Casimir element 
$\Delta \in \mathbf{C}[sl_2]$. Then we have a Poisson structure on $A$. 
Notice that $\mathrm{Spec}A$ is a 5-dimensional affine space $\mathbf{A}^5$.  
The Poisson centre $C_n := \{g \in A; \{g, A\} = 0\}$ is the polynomial ring 
$\mathbf{C}[f_n]$ generated by an element $f_n$ of $A$. 
The ring homomorphism $\mathbf{C}[f_n] \to A$ induces a morphism of 
algebraic varieties $f_n:  \mathbf{A}^5\to \mathbf{A}^1$. One can prove that 
$f_n$ coincides with the defining polynomial of $X_n$: $a^2x + 2aby + b^2z + (xz-y^2)^n$ 
after a suitable $\mathbf{C}^*$-equivariant coordinates change of 
$\mathbf{A}^5$. The Poisson structure on $\mathbf{A}^5$ induces a Poisson 
structure on the central fibre $X_n := \{f _n = 0\}$.  This Poisson structure has 
weight $-2$ and it is generically nondegenerate; in other words, $X_n$ admits 
a symplectic structure of weight $2$.  

Similarly, by using the symmetric product $S^3(\mathbf{C}^2)$ 
of the standard representation, we get a Poisson structure on $A$ 
with a Poisson centre $f$. Then $f$ is equivalent to the equation of the 
6-dimensional hypersurface.  

(c) The series $X_n$ of hypersurfaces can be also obtained as symplectic reductions of 
Hanany and Mekareeya [H-M] determined by unitrivalent graphs.    
 
Thus we have three symplectic structures on $X_n$ and have two symplectic structures 
on the 6-dimensional example.  
They are all equivalent by Theorem 3.1.  
\vspace{0.15cm}

In the latter half of the article we discuss the equivalence problem up to constant. 
Let $(X, \omega)$ be the same one as in Theorem (3.1); namely  
$l \ne 0$. 
A symplectic structure $\omega'$ on 
$X$ is equivalent to $\omega$ up to constant when $\omega' = \lambda\cdot \omega$ with 
some $\lambda \in \mathbf{C}^*$. If the weight $l$ of $\omega$ is nonzero, then 
the equivalence up to constant implies the equivalence up to $\mathbf{C}^*$-equivariant 
automorphism.  
Let $R$ be the affine ring of $X$. By the assumption $R$ is positively 
graded: $R = \oplus_{i \geq 0}R_i$. We put $\mathbf{P}(X) := \mathrm{Proj}(\oplus_{i \geq 0}R_i)$. 
Roughly speaking, we reduce the equivalence problem for the symplectic structure on $X$ 
to the uniqueness of the contact structure on $\mathbf{P}(X)$.
 
It is well known that a contact structure is an odd dimensional counterpart of a symplectic 
structure in complex and differential geometry. The author thinks that this is a 
good occasion to give an appropriate formulation of the contact structure for  
singular varieties. 

Recall that a contact structure on a complex manifold $Z$ of dimension $2d + 1$ is an 
exact sequence of vector bundles 
$$ 0 \to D \to TZ \stackrel{\theta}\to M \to 0, $$ 
with $\mathrm{rank}(D) = 2d$ and $\mathrm{rank}(M) = 1$ so that 
$d\theta\vert_D$ induces a non-degenerate pairing on $D$. 
The line bundle $M$ is called the contact line bundle. 
According to LeBrun [LeB], the contact structure is a unique one with the contact line 
bundle $M$ if and only if $H^0(Z, O(D)) = 0$. 
  
Let us consider the natural projection map $p: X - \{0\} \to \mathbf{P}(X)$. 
Then all fibres of $p$ are isomorphic to $\mathbf{C}^*$, but some of them are multiple fibres.   
There exists an open dense subset $\mathbf{P}(X)^0$ of $\mathbf{P}(X)$ such that $\mathbf{P}(X)^0$ 
is smooth and $p$ is a $\mathbf{C}^*$-bundle over $\mathbf{P}(X)^0$. Define  
$L := O_{\mathbf{P}(X)}(1)\vert_{\mathbf{P}(X)^0}$. 
The symplectic form $\omega$ on $X_{reg}$ of weight $l \ne 0$  
determines a contact structure on $\mathbf{P}(X)^0$ with the contact line bundle $L^{\otimes l}$ 
(cf. (4.3)).  
If $\mathrm{Codim}_{\mathbf{P}(X)}(\mathbf{P}(X) - \mathbf{P}(X)^0) \geq 2$, one can 
employ this contact structure on $\mathbf{P}(X)^0$ as a contact structure on $\mathbf{P}(X)$. 
But when $\mathrm{Codim}_{\mathbf{P}(X)}(\mathbf{P}(X) - \mathbf{P}(X)^0) = 1$, the contact structure on 
$\mathbf{P}(X)^0$  does not yet have enough information. This is the case, for example, when 
$(X, \omega)$ is a Du Val singularity with a symplectic structure of weight $2$, 
So, in a general case, 
we need to introduce the notion of a {\em contact orbifold structure} (see (4.4) for details). 
A contact orbifold 
structure on a normal variety $Z$ consists of an orbifold structure $Z^{orb}$ on $Z$, an orbifold 
line bundle $\mathcal{M}$ (= contact line bundle) on $Z^{orb}$ and a global section $\theta$ of 
$\underline{\mathrm{Hom}}(\Theta_{Z^{orb}}, \mathcal{M})$. Then one can prove:   
\vspace{0.15cm}

{\bf Theorem (4.4.1)}. {\em The projectivised cone $\mathbf{P}(X)$ has a contact orbifold structure.} 
\vspace{0.15cm}

Let $\mathcal{L} \in \mathrm{Pic}(\mathrm{P}(X)^{orb})$ be the tautologcal line bundle 
and assume that $\mathcal{M} = \mathcal{L}^{\otimes l}$. 
Then one can completely recover the original symplectic structure $(X, \omega)$
from the data $(\mathbf{P}(X), \mathcal{M}, \theta )$. 

For each Du Val singularity $(X, \omega)$ of type ADE, a contact orbifold structure on   
$\mathbf{P}(X) \cong \mathbf{P}^1$ is determined. 
But these structures are all different even though the underlying space is the same $\mathbf{P}^1$.  
In other words, $\mathbf{P}^1$ has infinitely many different contact orbifold structures. 

When $X$ has canonical singularities, the projectivised cone 
$\mathbf{P}(X)$ is a singular Fano variety. But $\mathbf{P}(X)$ turns out to be a 
very special one. In fact, we prove that the contact orbifold 
structure $(\mathbf{P}(X), \mathcal{M}, \theta )$ is rigid under a small  
deformation if $X$ has canonical singularities (Proposition (5.2)). 
When $X$ is the closure of a minimal nilpotent orbit $O_{min}$ of a simple Lie algebra, 
$\mathbf{P}(X)$ is a contact Fano homogeneous manifold. In this case the contact structure is known to be 
rigid under a small deformation (cf. [LeB]). Thus Proposition (5.2) generalises this fact.   

The equivalence problem for a symplectic structure on $X$ is now  
reduced to the uniqueness of the contact orbifold structure on 
$\mathbf{P}(X)$. In most examples the symplectic structures turn out to be unique up to 
contant with very few exceptions (Section 6).     
    
Section 7 is a speculation based on the analogy of the 
Bogomolov decomposition for compact K\"{a}hler manifolds with $c_1 = 0$. 
The contents of Section 6 are still fragmentary. However, the problems addressed in 
the final section would play a role as a working hypothesis in the future study.  

\vspace{0.2cm}

{\bf 1. Equivalence problem for complex analytic germs} 

Assume that $X$ is a germ of a normal complex space whose regular locus $X_{reg}$ 
admits a holomorphic symplectic 2-form $\omega$.  
Two such pairs $(X,\omega)$ and $(X', \omega')$ 
are {\em equivalent} if there is an isomorphism $\phi : X \to X'$ such that 
$\omega = \phi^*(\omega')$.  They are not, a priori, equivalent 
even if their underlying complex analytic structure are equivalent. The Darboux theorem  
asserts that any holomorphic symplectic structure on $(\mathbf{C}^{2n},0)$ is 
equivalent to the standard one $dx_1 \wedge dx_2 + ... + dx_{2n-1} \wedge dx_{2n}$.

One can generalise the Darboux theorem to a quotient singularity, which might 
be already known. 
\vspace{0.2cm}

{\bf Proposition (1.1)}.  {\em Let $(X,0)$ be a quotient symplectic singularity with a holomorphic symplectic 
form $\omega$. Then any holomorphic symplectic form on $(X,0)$ is equivalent to $\omega$. } 
          
{\em Proof}.  Write $X = \mathbf{C}^{2n}/G$ with a finite group $G \subset Sp(2n, \mathbf{C})$.  
Let $\pi: (\mathbf{C}^{2n},0) \to (X,0)$ 
be a natural projection. Let $\omega'$ be an arbitrary symplectic form on $(X,0)$. Let $\tilde{\omega}$ 
and $\tilde{\omega'}$ be respectively the pull-backs of $\omega$ and $\omega'$ by $\pi$. 
We shall prove that there is a $G$-equivariant isomorphism $\tilde{\varphi} : (\mathbf{C}^{2n},0) \to 
(\mathbf{C}^{2n},0)$ such that $\tilde{\varphi}^*(\tilde{\omega'}) = \tilde{\omega}$.  
Then this $\tilde{\varphi}$ descends to an automorphism $\varphi$ of $(X,0)$ such that 
${\varphi}^*(\omega') = \omega$. 
We first prove a linear algebra version of this fact: 
\vspace{0.2cm}

{\bf  Lemma (1.2)}. {\em Let $V$ be a $2n$-dimensional complex representation of a finite 
group $G$. Assume that $\omega$ and $\omega'$ are  $G$-invariant non-degenerate 
skew-symmetric 2-forms on $V$. Then there is a $G$-equivariant linear isomorphism 
$\phi$ such that $\phi^*(\omega) = \omega'$. }
\vspace{0.15cm}

{\em Proof}. 
Denote by $V^*$ the dual representation of $V$.  
We divide irreducible representations $V$ of $G$ into three types: 

(I)  $V \cong V^*$ and  $\dim (\wedge^2 V^*)^G = 1$,   

(II) $V \cong V^*$ and $(\wedge^2 V^*)^G = 0$, 

(III) $V$ is not isomorphic to $V^*$ as a $G$-module. 

Note that if $V$ is irreducible and $V \cong V^*$, then $\mathbf{C} = \mathrm{Hom}_G(V,V) = (V \otimes V^*)^G 
 = (V^* \otimes V^*)^G = (\wedge^2 V^*)^G \oplus (Sym^2(V^*))^G$.  In case (I) one has 
$\dim (\wedge^2 V^*)^G = 1$ and the isomorphism $V \cong V^*$ is given by a 
$G$-invariant non-degenerate skew-symmetric form which is unique up to scalar. 
In case (II) one has $\dim (Sym^2(V^*))^G = 1$ and $V \cong V^*$ is given by a $G$-invariant 
non-degenerate symmetric form which is unique up to scalar. 
If $V$ is of type (III), then $(\wedge^2 V^*)^G = 0$ because there is 
an injection $(\wedge^2 V^*)^G \to (V^* \otimes V^*)^G = \mathrm{Hom}_G(V, V^*) = 0$. 
Moreover $\dim (V \otimes V^*)^G = 1$ because $(V \otimes V^*)^G = \mathrm{Hom}_G(V,V) 
= \mathbf{C}$. 
Finally note that if $V$ and $V'$ are irreducible representations of different type, 
one has $(V \otimes V')^G =0$ and $\mathrm{Hom}_G(V, V') = 0$. 
         
Assume that $V$ is of type I. An element  
$$\varphi \in (\wedge^2V^*)^G = \mathrm{Hom}_G(V, V^*)$$ is 
represented by a matrix $X$ if we choose a basis of $V$ and choose its 
dual basis of $V^*$.    
By changing the initial basis if necessarily, we may assume that 
$X = aJ$, where 
$a \in \mathbf{C}$ and where   
$$ J = \left( \begin{array}{cc}
0 & I  \\ 
-I & 0  
\end{array}\right). $$ Here notice that $V$ is an even dimensional 
$\mathbf{C}$-vector space.   
Similarly an element $$\varphi \in (\wedge^2 ((V^*)^{\oplus n}))^G 
\subset \mathrm{Hom}_G(V^{\oplus n}, (V^*)^{\oplus n})$$ 
is represented by a matrix 
$$  X = \left( \begin{array}{cccc}
a_{11}J & a_{12}J & \cdots & a_{1n}J \\ 
a_{21}J & a_{22}J & \cdots & a_{2n}J \\ 
\cdots & \cdots & \cdots & \cdots \\
a_{n1}J & a_{n2}J & \cdots & a_{nn}J 
\end{array}\right),$$ 
where $A := (a_{ij})$ is a symmetric matrix. 
If $\varphi$ is nondegenerate, then for a suitable matrix $T$ of the form 
$$ T = \left( \begin{array}{cccc}
t_{11}I & t_{12}I & \cdots & t_{1n}I \\ 
t_{21}I & t_{22}I & \cdots & t_{2n}I \\ 
\cdots & \cdots & \cdots & \cdots \\
t_{n1}I & t_{n2}I & \cdots & t_{nn}I 
\end{array}\right),$$ 
we have 
$$ {}^tTXA = \left( \begin{array}{cccc}
J & 0& \cdots & 0\\ 
0 & J & \cdots & 0 \\ 
\cdots & \cdots & \cdots & \cdots \\
0 & 0 & \cdots & J 
\end{array}\right).$$ 
 
Assume that $V$ is of type II. 
Then an element  
$$\varphi \in (Sym^2V^*)^G = \mathrm{Hom}_G(V, V^*)$$ can be  
represented by a matrix of the form $aI$ with $a \in \mathbf{C}^*$. 
Similarly an element $$\varphi \in (Sym^2 ((V^*)^{\oplus n}))^G 
\subset \mathrm{Hom}_G(V^{\oplus n}, (V^*)^{\oplus n})$$ 
is represented by a matrix 
$$  X = \left( \begin{array}{cccc}
a_{11}I & a_{12}I & \cdots & a_{1n}I \\ 
a_{21}I &  & \cdots & a_{2n}I \\ 
\cdots & \cdots & \cdots & \cdots \\
a_{n1}I & a_{n2}I & \cdots & a_{nn}I 
\end{array}\right),$$ 
where $A := (a_{ij})$ is a skew-symmetric matrix. 
If $\varphi$ is nondegenerate, then for a suitable matrix $T$ of the form 
$$ T = \left( \begin{array}{cccc}
t_{11}I & t_{12}I & \cdots & t_{1n}I \\ 
t_{21}I & t_{22}I & \cdots & t_{2n}I \\ 
\cdots & \cdots & \cdots & \cdots \\
t_{n1}I & t_{n2}I & \cdots & t_{nn}I 
\end{array}\right),$$ 
we have 
$$ {}^tTXA = \left( \begin{array}{cc}
0 & I\\ 
-I & 0 
\end{array}\right).$$ 
Finally assume that $V$ is of type III. 
Then $(\wedge^2(V \oplus V^*))^G \subset \mathrm{Hom}_G(V \oplus V^*, V^* \oplus V)$. 
An element $(\wedge^2(V \oplus V^*))^G$ is represented by a matrix 
$$ \left( \begin{array}{cc}
0 & aI\\ 
-aI & 0 
\end{array}\right).$$ 
Similarly an element 
$$\varphi \in (\wedge^2(V^{\oplus n} \oplus (V^*)^{\oplus n}))^G 
\subset \mathrm{Hom}_G(V^{\oplus n} \oplus (V^*)^{\oplus n}, 
(V^*)^{\oplus n} \oplus V^{\oplus n})$$ is represented by a matrix 
$$  X = \left( \begin{array}{cc}
0 & A\\ 
-{}^tA &  0  
\end{array}\right), $$ 
where 
$$  A = \left( \begin{array}{cccc}
a_{11}I & a_{12}I & \cdots & a_{1n}I \\ 
a_{21}I &  & \cdots & a_{2n}I \\ 
\cdots & \cdots & \cdots & \cdots \\
a_{n1}I & a_{n2}I & \cdots & a_{nn}I 
\end{array}\right).$$ 
If $\varphi$ is nondegenerate, then for a suitable matrix $T$ of the form 
$$ \left( \begin{array}{cc}
T_1 & 0 \\ 
0 & T_1  
\end{array}\right),$$ 
we have 
$$ {}^tTXA = \left( \begin{array}{cc}
0 & I\\ 
-I & 0 
\end{array}\right).$$ 

Now let us consider the $V$ in Lemma.  Decompose $V$ into the sum of irreducible 
representations 
$$ V = \bigoplus (V_i)^{\oplus l_i} \oplus \bigoplus (V'_j)^{\oplus m_j} 
 \oplus \bigoplus (W_k)^{\oplus n_k}, $$      
where $V_i$ are of type (I), $V'_j$ are of type (II) and $W_k$ are of type (III). 
Since $V$ admits a $G$-invariant non-degenerate 2-form $\varphi$, we see that, in 
the third factor $\bigoplus (W_k)^{\oplus n_k}$ each irreducible representation 
and its dual one appear in a pairwise way. Thus the third factor can be written 
as $\bigoplus (W_k \oplus W_k^*)^{\oplus n_k}$.  
By the observations above we see that $\varphi$ can be transformed to a   
standard $G$-equivariant 2-form after making a suitable $G$-equivariant base 
change of $V$.   Q.E.D.   
\vspace{0.2cm}

Let us return to the proof of Proposition (1.1).  Let $\tilde{\omega}(0) \in \wedge^2T^*_0({\mathbf C}^{2n})$ 
and $\tilde{\omega'}(0) \in \wedge^2T^*_0({\mathbf C}^{2n})$ be respectively the restriction of 
$\tilde{\omega}$ and $\tilde{\omega'}$ to the origin $0 \in \mathbf{C}^{2n}$. 
By the lemma above, we may assume from the first that $\tilde{\omega}(0) = 
\tilde{\omega'}(0)$. The rest of the argument is an equivariant version of Moser's  
standard argument.
For $\tau \in \mathbf{R}$, define  
$$ \tilde{\omega}_{\tau} :=  (1-\tau)\tilde{\omega} + \tau\tilde{\omega'}. $$
We put $$ u := d\tilde{\omega}_{\tau}/d\tau.$$ 
Let us consider the complex 
$(\pi^G_*\Omega^{\cdot}_{\mathbf{C}^{2n}}, d)$, which is a resolution of the consant sheaf 
$\mathbf{C}_X$. 
Note that $u$ is a section of $\pi^G_*\Omega^2_{\mathbf{C}^{2n}}$. 
Since $u$ is d-closed, one can write $u = dv$ with a $G$-invariant 1-form $v$. Moreover, 
$v$ can be chosen such that $v(\mathbf{0}) = 0$. 
Define a vector field $X_{\tau}$ on $(\mathbf{C}^{2n},0)$ by 
$$ i_{X_{\tau}} \tilde{\omega}_{\tau} = -v. $$
Since $\tilde{\omega}(\tau)$ is $d$-closed, we have 
$$ L_{X_{\tau}}\tilde{\omega}_{\tau} = -u, $$ 
where $L_{X_{\tau}}\tilde{\omega}_{\tau}$ is the Lie derivative of $\tilde{\omega}_{\tau}$ along 
$X_{\tau}$. If we take a sufficiently small open set $V$ of $\mathbf{0} \in \mathbf{C}^{2n}$, 
then the vector fields $\{X_{\tau}\}_{0 \le \tau \le 1}$ define a family of open 
immersions $\tilde{\varphi}_{\tau}: V \to \mathbf{C}^{2n}$ via 
$$ d\tilde{\varphi}_{\tau}/d\tau = X_{\tau}(\tilde{\varphi}_{\tau}), \; \tilde{\varphi}_0 = id.$$ 
Since all $\tilde{\varphi}_{\tau}$ fix the origin and $X_{\tau}$ are all $G$-invariant, 
$\tilde{\varphi}_{\tau}$ induce $G$-invariant automorphisms of $(\mathbf{C}^{2n},0)$. 
We have 
$$ d(\tilde{\varphi}_{\tau}^*\tilde{\omega}_{\tau})/d\tau = \tilde{\varphi}_{\tau}^*(d\tilde{\omega}_{\tau}/d\tau 
+ L_{X_{\tau}}\tilde{\omega}_{\tau}) = 0.$$ 
In particular, $\tilde{\varphi}_0^*\tilde{\omega}_0 = \tilde{\varphi}_1^*\tilde{\omega}_1$. 
The left hand side is $\tilde{\omega}$ and right hand side is $\tilde{\varphi}_1^*\tilde{\omega'}$.  
If we put $\tilde{\varphi} := \tilde{\varphi}_1$, then $\tilde{\varphi}$ is a desired 
$G$-equivariant automorphism of $(\mathbf{C}^{2n},0)$. Q.E.D.         
\vspace{0.2cm}

{\bf 2. Affine varieties with $\mathbf{C}^*$-actions and symplectic structures} 

Let $X$ be a normal affine variety of 
dimension $2d$ with a $\mathbf{C}^*$-action. Assume that $0 \in X$ is a unique fixed point of the $\mathbf{C}^*$-action with positive weights. More precisely, the cotangent space $m_{X,0}/m^2_{X,0}$ of $0 \in X$ has 
only positive weights with respect to the $\mathbf{C}^*$-action or equivalently, the affine ring $R$ of $X$ is positively graded: $\oplus_{i \ge 0}R_i$ with $R_0 = \mathbf{C}$. 
In the remainder we call such a $\mathbf{C}^*$-action a {\em good} $\mathbf{C}^*$-action.  
Let $\omega$ be an algebraic symplectic 
2-form on $X_{reg}$ with weight $l$. If we represent the $\mathbf{C}^*$-action by 
the family $\{\phi_t\}_{t \in \mathbf{C}^*}$ of automorphisms of $X$, then $\phi_t^*(\omega) = t^l \cdot \omega$.  

\vspace{0.2cm}

{\bf Lemma (2.1)}. {\em If $\omega'$ is another symplectic 2-form with 
 weight $l'$, then $l = l'$.}
\vspace{0.2cm}

{\em Proof}.  Assume that $l < l'$. Since ${\omega'}^d$ is 
a generator of the canonical line bundle $K_X$, one can write $\omega^d = g\cdot {\omega'}^d$ 
with a homogeneous regular function $g$ on $X$ with negative weight $l-l'$. But this contradicts 
the assumption that $X$ is positively weighted.
\vspace{0.2cm}

{\bf Remark}. The lemma shows that if we fix a $\mathbf{C}^*$-action on $X$, then 
$l$ is uniquely determined. But if we replace the $\mathbf{C}^*$-action on $X$ by a 
different one, $l$ may possibly change. For example, let $X$ be a 2-dimensional 
quotient singularity $\mathbf{C}^2/G$ where $G$ is a cyclic group of order $m$ 
acting on $\mathbf{C}^2$ as $x \to \zeta \cdot x$ and $y \to \zeta^{-1}\cdot y$ 
with a primitive $m$-th root $\zeta$ of unity. Introduce a $\mathbf{C}^*$-action 
on $\mathbf{C}^2$ by $x \to t^p \cdot x$ and $y \to t^q \cdot y$ with positive 
integers $p$ and $q$ which are coprime to each other. Put $u := x^m$, $v := y^m$ 
and $w := xy$. Then $X$ is an affine subvariety of $\mathbf{C}^3(u,v,w)$ defined by the equation 
$uv - w^m = 0$. The $\mathbf{C}^*$-action on $\mathbf{C}^2(x,y)$ descends to a 
$\mathbf{C}^*$-action on $X$. With respect to this $\mathbf{C}^*$-action, we have 
$$(wt(u), wt(v), wt(w)) = (mp, mq, p+q).$$ 
If we choose $p$, $q$ in such a way that $p+q$ and $m$ are coprime, then 
$GCD (mp, mq, p+q) = 1$. By definition $X$ has a symplectic 2-form 
$$ \omega := du \wedge dv/ w^{m-1}, $$ which has weight $p + q$. 
\vspace{0.2cm}     

Before going to the next lemma, we recall the notions of a symplectic singularity and a 
canonical singularity. Let $(X, \omega)$ be a normal affine variety with a $\mathbf{C}^*$-action and an algebraic symplectic 2-form $\omega$ with weight $l$.  Since $\omega^d := \omega \wedge ... \wedge \omega$ is a generator of the dualizing sheaf $\omega_X$, the canonical divisor $K_X$ is a Cartier divisor. Let $\pi: Y \to X$ be a resolution and let $E_i$ ($1 \le i \le n$) be the $\pi$-exceptional divisors. 
One can write $K_Y = \pi^*K_X + \Sigma a_iE_i$ with 
some integers $a_i$. If $a_i \geq 0$ for all $i$, then we say that $X$ has {\em canonical 
singularities}. On the other hand, if $\omega$ is pulled back to a regular 2-form on 
$Y$, we say that $X$ has {\em symplectic singularities} [Be]. By [Na 2]  $X$ has canonical singularities if and only if 
$X$ has symplectic singularities.

In order to check that $X$ does not have canonical singularities, we only have 
to find a {\em partial} resolution $f: Z \to X$ such that $\mathrm{Exc}(f)$ contains a 
divisor $E$ such that $f^*(\omega^d)$ has a pole along $E$.    
\vspace{0.2cm}

{\bf Lemma (2.2)}. {\em If $X$ has only canonical singularities, then $l$ is positive.} 
\vspace{0.2cm}

{\em Proof}. We prove that if $l \leq 0$, then $X$ does not have canonical singularities. 
Let $R$ be the affine ring of $X$. By the $\mathbf{C}^*$-action of $X$, $R$ has a 
grading $R = \oplus_{k \geq 0}R_k$ with $R_0 = \mathbf{C}$. Let $x_0$, ..., $x_n$ be 
homogeneous minimal generators of the $\mathbf{C}$-algebra $R$ and put $a_i := wt(x_i)$. 
We assume that $\mathrm{GCD}(a_0, ..., a_n) = 1$. The affine variety $X$ is embedded in 
$\mathbf{C}^{n+1}$ by $x_i$'s. Let $\pi: V \to \mathbf{C}^{n+1}$ be the weighted blowing 
up of $\mathbf{C}^{n+1}$ with weight $(a_0, ..., a_n)$. By the definition, 
$V$ is covered by open sets $V_i$ $(0 \le i \le n)$ and there is a $\mathbf{Z}/a_i\mathbf{Z}$-Galois 
cover 
$$ p_i: \mathbf{C}^{n+1} \to V_i$$ 
such that 
$$ (\pi \circ p_i)^*x_i = (x'_i)^{a_i} $$ 
$$ (\pi \circ p_i)^*x_j = (x'_i)^{a_j}x'_j \; (j \ne i), $$ 
and $p_i$ is the quotient map of the $\mathbf{Z}/a_i\mathbf{Z}$-action on 
$\mathbf{C}^{n+1}$ 
$$ x'_i \to \zeta \cdot x'_i, $$ 
$$ x'_j \to \zeta^{-a_j} \cdot x'_j $$ 
with an $a_i$-th primitive root $\zeta$ of unity. 
The exceptional divisor $E := \pi^{-1}(0)$ is isomorphic to 
the weighted projective space $\mathbf{P}(a_0, ..., a_n)$. 
Let us observe the restriction of $p_i$ to $p_i^{-1}(E \cap V_i)$. 
Note that $p_i^{-1}(E \cap V_i)$ is a divisor of $\mathbf{C}^{n+1}$ defined 
by the equation $x'_i = 0$ and that the $\mathbf{Z}/a_i\mathbf{Z}$-action on 
$p_i^{-1}(E \cap V_i)$ is given by 
$$ x'_j \to \zeta^{-a_j}\cdot x'_j.$$ 
By the assumption $GCD(a_0, ..., a_n) = 1$, we see that $\mathbf{Z}/a_i\mathbf{Z}$ 
acts effectively on $p_i^{-1}(E \cap V_i)$. Therefore, 
$$ p_i^{-1}(E \cap V_i) \to E \cap V_i $$ 
is a $\mathbf{Z}/a_i\mathbf{Z}$-Galois covering. Let $p \in E$ be a general point. 
Then $V$ is smooth at $p$. Let $\tilde{X}$ be the proper transform 
of $X \subset \mathbf{C}^{n+1}$ by the weighted blowing up $\pi: V \to \mathbf{C}^{n+1}$ 
and let 
$$ \pi_X: \tilde{X} \to X$$ 
be the induced birational morphism. 
Note that 
$$ E \cap \tilde{X} = \mathrm{Proj}(\oplus_{k \geq 0}R_k).$$ 
Since $E \cap \tilde{X}$ is generically smooth and $E$ is a Cartier divisor at a 
general point $p \in E \cap \tilde{X}$, we can see that $\tilde{X}$ is also smooth 
at such a point $p$. 

Now let us consider the $2d$-form $\omega^d$ and regard it as a section 
of the canonical line bundle $K_X$. 
We shall prove that $(\pi_X)^*\omega^d$ has a pole along $E \cap \tilde{X}$ if $l \leq 0$. 
Take a general point $p \in E \cap \tilde{X}$ and assume that $p \in V_i$. 
We put $\tilde{X}_i := (p_i)^{-1}(\tilde{X} \cap V_i)$ and $E_i := (p_i)^{-1}(E \cap V_i)$. 
Recall that $p_i^{-1}(E \cap V_i) \to E \cap V_i$ is a $\mathbf{Z}/a_i\mathbf{Z}$-Galois 
covering whose branch locus is contained in the divisor $\prod_{j \ne i}x_j = 0$ 
of $E = \mathbf{P}(a_0, ..., a_n)$. 
Since $\mathrm{Proj}(\oplus_{k \geq 0}R_k)$ is not contained in the divisor 
$\prod x_j = 0$ of $\mathbf{P}(a_0, ..., a_n)$, we see that 
$$ E_i \cap \tilde{X}_i \to E \cap \tilde{X} \cap V_i$$ 
is a $\mathbf{Z}/a_i\mathbf{Z}$-Galois cover. 
This implies that the order of the zeros (or the poles) of $(\pi_X)^*\omega^d$ 
along $E \cap \tilde{X}$ coincides with the order of the zeros (or the poles) of 
$(\pi_X \circ p_i\vert_{{\tilde X}_i})^*\omega^d$ along 
$E_i \cap \tilde{X}_i$. Let $q \in \tilde{X}_i$ be a point such that 
$p_i(q) = p$. One can choose the local coordinates of $q \in \tilde{X}_i$ 
from $x'_j - x'_j(q)$ $(0 \le j \le n)$. Since $E_i$ is smooth at $q$, 
we can include $x'_i$ among the local coordinates (note that $x'_i(q) = 0$). 
Assume that $x'_i$, $x'_{j_1} - x'_{j_1}(q)$, ..., 
$x'_{j_{2d-1}} - x'_{j_{2d-1}}(q)$ are local coordinates. 
Recall that $V$ has a natural $\mathbf{C}^*$-action and this $\mathbf{C}^*$-action 
extends to the $\mathbf{C}^*$-action on $(x'_0, ..., x'_n) \in \mathbf{C}^{n+1}$ 
by $$ x'_i \to t \cdot x'_i, $$ 
and $$ x'_j \to x'_j, \; (j \ne i).$$ 
Since $\omega$ has weight $l$, the weight of $(\pi_X\circ p_i\vert_{\tilde{X}_i})^*\omega^d$ 
is $d \cdot l$. Around $q \in \tilde{X}_i$, one can write 
$$ (\pi_X \circ p_i\vert_{\tilde{X}_i})^*\omega^d = h \cdot dx'_i \wedge 
dx'_{j_1} \wedge ... \wedge dx'_{j_{2d-1}}$$ 
with a meromorphic function $h$ of degree $d \cdot l - 1$. This means that 
$(\pi_X \circ p_i\vert_{\tilde{X}_i})^*\omega^l$ has poles of order $1 - d \cdot l$ along 
$E_i \cap \tilde{X}_i$ if $l \leq 0$. Q.E.D. 
\vspace{0.2cm}

{\bf Corollary (2.3)}. {\em If $\mathrm{Codim}_X \mathrm{Sing}(X) \geq 4$, then $l > 0$.} 
\vspace{0.15cm}

{\em Proof}. If $\mathrm{Sing}(X)$ has at least codimension $4$ in $X$, then the symplectic 2-form $\omega$ 
extends to a regular 2-form on 
an arbitrary resolution $\tilde{X}$ of $X$ by Flenner [Fl]. This implies that $X$ has only canonical 
singularities. 
\vspace{0.15cm}

{\bf Lemma (2.4)}. {\em If $X$ has an isolated singularity and $l > 0$, then $X$ has only canonical 
singularities.} 
\vspace{0.2cm}

{\em Proof}. Let $\pi: Y \to X$ be a $\mathbf{C}^*$-equivariant resolution. 
Let $Y_c$ be a relatively compact open 
subset of $Y$ such that $\pi^{-1}(0) \subset Y_c$. 
Write $K_Y = \pi^*K_X + \Sigma a_iE_i$ where $E_i$ are $\pi$-exceptional 
divisors. Since $K_X$ is Cartier (because of the existence of $\omega$), 
all coefficients $a_i$ are integers. In order to prove that $a_i \geq 0$, 
we only need to prove that $a_i > -1$.  This condition is equivalent to the 
$L^2$-condition (cf. the proof of [Ko, Proposition 3.20]):   
$$ \int_{Y_c} \pi^*\omega^d \wedge \pi^*{\bar{\omega}}^d < \infty.$$ 
Since $\mathbf{R}_{>0}$ is naturally contained in 
$\mathbf{C}^*$, each element $t \in \mathbf{R}_{>0}$ acts on $X$ as an 
automorphism $\phi_t$ of $X$.   
Let $U$ be an open neighborhood of $0 \in X$ such that $\phi_t(U) \subset 
U$ for all $t \in (0, \; 1]$. Put $V := \pi^{-1}(U)$. Fix $\epsilon_0 \in (0, \; 1)$ and put 
$U_n := \phi_{\epsilon_0^n}(U) - \phi_{\epsilon_0^{n+1}}(U)$. 
Define $V_n := \pi^{-1}(U_n)$.        
Since $\phi_t^*\omega = t^l \cdot \omega$, we have 
$$ \int_{V_n} \pi^*\omega^d \wedge \pi^*{\bar{\omega}}^d = 
\epsilon_0^{2dnl} \cdot \int_{V_0} \pi^*\omega^d \wedge \pi^*{\bar{\omega}}^d.$$ 
By the definition we have 
$$ \int_V \pi^*\omega^d \wedge \pi^*{\bar{\omega}}^d 
= \sum_{n=0}^{\infty} \int_{V_n} \pi^*\omega^d \wedge \pi^*{\bar{\omega}}^d.$$ 
But the right hand side equals $$(\sum_{n=0}^{\infty}\epsilon_0^{2dnl})\int_{V_0}
\pi^*\omega^d \wedge \pi^*{\bar{\omega}}^d < \infty.$$ 
The desired $L^2$-condition has now been proved. 

Note that this proof is not valid for a non-isolated case because 
$ \int_{V_n} \pi^*\omega^d \wedge \pi^*{\bar{\omega}}^d$ might be 
infinite.  Q.E.D.             
\vspace{0.2cm}

{\bf 3. Algebraic version of equivalence problems} 

In this section $(X, \omega)$ is a pair of a normal affine variety $X$ of dimension $2d$ with a good $\mathbf{C}^*$-action and an algebraic symplectic 2-form $\omega$ 
on $X_{reg}$ with weight $l$.  
We shall consider the equivalence problem for a pair 
$(X, \omega)$. Let $(X', \omega')$ be another pair. Then $(X, \omega)$ and $(X', \omega')$ 
are equivalent if there is a $\mathbf{C}^*$-equivariant isomorphism $\phi: X \cong X'$ such that $\omega = \phi^*(\omega')$.  In particular, if $X = X'$, then $\omega$ and $\omega'$ are 
called equivalent symplectic structures on $X$. A purpose of this section is to prove 
\vspace{0.2cm}

{\bf Theorem (3.1)}. {\em Assume that $l \ne 0$. Then  $\omega$ is a unique  
symplectic structure with weight $l$ on $X$ up to equivalence. }      
\vspace{0.2cm}

We shall briefly recall some basic results on Poisson structures and their deformations. 
For details see [Na 1].  
Note that the symplectic 2-form $\omega$ gives a natural Poisson structure $\{\;, \;\}$ on 
$X_{reg}$. By the normality of $X$, this Poisson structure extends to a Poisson 
structure $X$. We denote this bracket also by $\{\;, \;\}$. 
The bracket $\{\;, \;\}$ has weight $-l$ with respect to the $\mathbf{C}^*$-action because 
$\omega$ has weight $l$. 
Namely if $f$ and $g$ are homogeneous element of $O_X$ of degree $a$ and $b$, then 
$\{f, g\}$ is a homogeneous element of degree $a + b - l$.

By using the Poisson bracket we define the Lichnerowicz-Poisson complex 
$$ 0 \to \Theta_{X_{reg}} \stackrel{\delta_1}\to \wedge^2 \Theta_{X_{reg}} \stackrel{\delta_2}\to  
... $$ 
by 
$$ \delta_p f(da_1 \wedge ... \wedge da_{p+1})  
:= \sum_{i = 1}^{p+1} (-1)^{i+1}\{a_i, f(da_1 \wedge 
... \hat{da_i} \wedge ... \wedge da_{p+1}\} $$   
$$ + \sum_{j < k} (-1)^{j+k}f(d\{a_j, a_k\} \wedge da_1 \wedge ... \wedge 
\hat{da_j}  
\wedge ... \wedge \hat{da_k} \wedge ... \wedge da_{p+1}).$$ 
In the Lichnerowicz-Poisson complex, $\wedge^p \Theta_{X_{reg}}$ is placed in degree $p$. 
By the symplectic form $\omega$ each term $\wedge^p \Theta_{X_{reg}}$ 
can be identified with the sheaf $\Omega^p_{X_{reg}}$ of $p$-forms. 
Moreover the Lichnerowicz-Poisson complex is identified with the truncated 
De Rham complex 
$$ 0 \to \Omega^1_{X_{reg}} \stackrel{d}\to \Omega^2_{X_{reg}} \stackrel{d}\to ... $$   

Put $S_1 := \mathrm{Spec}\mathbf{C}[\epsilon]$. Then the 2-nd cohomology 
$\mathbf{H}^2(\Gamma(X_{reg}, \wedge^{\ge 1}\Theta_{X_{reg}}))$ describes the 
equivalence classes of the $O_{S_1}$-bilinear Poisson structures 
$\{\;, \;\}_{\epsilon}$ on $X_{reg} \times S_1$ which are extensions of the original Poisson structure 
$\{\;, \;\}$ on $X_{reg} \times \{0\}$.  
In fact, for $\varphi \in \Gamma (X_{reg}, \wedge^2 \Theta_{X_{reg}})$, we define 
a bracket $\{\;, \;\}_{\epsilon}$ on $O_{X_{reg}} \oplus \epsilon O_{X_{reg}}$ 
by $$ \{f + \epsilon f', g + \epsilon g'\}_{\epsilon} := \{f,g\} + \epsilon (\varphi(df \wedge dg) + 
\{f, g'\} + \{f', g\}).$$ 
Then this bracket is a Poisson bracket if and only if $\delta(\varphi) = 0$. 
On the other hand, an element $\theta \in \Gamma (X_{reg}, \Theta_{X_{reg}})$ corresponds 
to an automorphism $\phi_{\theta}$ of $X_{reg} \times S_1$ over $S_1$ which restricts to give 
the identity map of $X_{reg} \times \{0\}$. Let $\{\; , \;\}_{\epsilon, 1}$ and 
$\{\; , \;\}_{\epsilon,2}$ be the Poisson structures determined respectively by elements $\varphi_1$ and 
$\varphi_2$ of $\Gamma (X_{reg}, \wedge^2 \Theta_{X_{reg}})$. Then the two Poisson structures are 
equivalent under $\phi_{\theta}$ if $\varphi_1 - \varphi_2 = \delta(\theta)$. 

Note that a Poisson structure $\{\;, \;\}_{\epsilon}$ on $X_{reg} \times S_1$ uniquely 
extends to a Poisson structure on $X \times S_1$. This means that 
$\mathbf{H}^2(\Gamma(X_{reg}, \wedge^{\ge 1}\Theta_{X_{reg}}))$ also describes 
equivalence classes of the $O_{S_1}$-bilinear Poisson structures 
$\{\;, \;\}_{\epsilon}$ on $X \times S_1$ which are extensions of the original Poisson structure 
$\{\;, \;\}$ on $X \times \{0\}$.
 
Let us introduce a $\mathbf{C}^*$ action on $X \times S_1$ in such a way that 
it acts on the first factor by the original action and acts trivially on the second factor.    
The following proposition is a $\mathbf{C}^*$-equivariant version of the above observation.  
\vspace{0.2cm}  

{\bf Proposition (3.2)} (Rigidity proposition). {\em  
Let $\{\; , \;\}_{\epsilon, 1}$ and $\{\; , \;\}_{\epsilon, 2}$ be two Poisson structures on 
$X \times S_1$ relative to $S_1$, both of which have weight $-l \ne 0$ and induce the original Poisson structure on $X \times \{0\}$. Then there is a $\mathbf{C}^*$-equivariant automorphism of $X \times S_1$ over $S_1$ such that it induces the identity map of $X \times \{0\}$ and it sends $\{\; , \;\}_{\epsilon, 1}$ to $\{\; , \;\}_{\epsilon, 2}$.}
\vspace{0.2cm}

{\em Proof}. Let $(\wedge^{\geq 1}\Theta_{X_{reg}}, \delta)$ be the Lichnerowicz-Poisson complex for 
a Poisson manifold $X_{reg}$. The algebraic torus $\mathbf{C}^*$ acts on 
$\Gamma(X_{reg}, \wedge^p\Theta_{X_{reg}})$ and there is an associated grading 
$$ \Gamma(X_{reg}, \wedge^p\Theta_{X_{reg}}) = \oplus_{n \in \mathbf{Z}} \;  
\Gamma(X_{reg}, \wedge^p\Theta_{X_{reg}})(n).$$ 

The coboudary map $\delta$ has degree $-l$; thus we have a complex 
$$ \Gamma(X_{reg}, \wedge^1\Theta_{X_{reg}})(0) \stackrel{\delta_1}\to 
\Gamma(X_{reg}, \wedge^2\Theta_{X_{reg}})(-l) \stackrel{\delta_2}\to 
\Gamma(X_{reg}, \wedge^3\Theta_{X_{reg}})(-2l).$$ 
The middle cohomology $\mathrm{Ker}(\delta_2)/\mathrm{Im}(\delta_1)$ of this complex describes the equivalence classes  
of the extension of the Poisson structure $\{\; , \;\}$ on $X_{reg}$ to that on $X_{reg} \times S_1$ 
with weight $-l$ up to $\mathbf{C}^*$-equivariant automorphism of $X_{reg} \times S_1$ 
over $S_1$ that induces the identity map of $X_{reg} \times \{0\}$. 
Since each Poisson   
structure $X_{reg} \times S_1$ uniquely extends to that on $X \times S_1$,  
$\mathrm{Ker}(\delta_2)/\mathrm{Im}(\delta_1)$ also describes 
the equivalence classes of the  
extension of the Poisson structure $\{\; , \;\}$ on $X$ to that on $X \times S_1$ 
with weight $-l$ up to $\mathbf{C}^*$-equivariant automorphism of $X \times S_1$ 
over $S_1$ that induces the identity map of $X \times \{0\}$. 

The Lichnerowicz-Poisson complex $(\wedge^{\geq 1}\Theta_{X_{reg}}, \delta)$ is identified 
with the truncated De Rham complex $(\Omega_{X_{reg}}^{\geq 1}, d)$ by the symplectic 
form $\omega$.  
The algebraic torus $\mathbf{C}^*$ acts on 
$\Gamma(X_{reg}, \Omega^p_{X_{reg}})$ and there is an associated grading 
$$ \Gamma(X_{reg}, \Omega^p_{X_{reg}}) = \oplus_{n \in \mathbf{Z}} \; 
\Gamma(X_{reg}, \Omega^p_{X_{reg}})(n).$$
The coboundary map $d$ has degree $0$; thus we have a complex 
$$ \Gamma(X_{reg}, \Omega^1_{X_{reg}})(l) \stackrel{d_1}\to 
\Gamma(X_{reg}, \Omega^2_{X_{reg}})(l) \stackrel{d_2}\to 
\Gamma(X_{reg}, \Omega^3_{X_{reg}})(l).$$ 
Since $\omega$ has weight $l$, this complex is identified with 
the 3 term complex above. 

We shall prove that $\mathrm{Ker}(d_2)/\mathrm{Im}(d_1) = 0$.     
The $\mathbf{C}^*$-action on $X$ defines a vector field 
$\zeta$ on $X_{reg}$. According to Naruki [Naru, Lemma 2.1.1] we define   
$$ \Delta: \Gamma(X_{reg}, \Omega^2_{X_{reg}}) \to \Gamma (X_{reg}, \Omega^1_{X_{reg}}) $$ 
by $\Delta (v) := i_{\zeta} v$. Since $\zeta$ is a $\mathbf{C}^*$-equivariant vector field, 
$\Delta$ induces a map 
$$ \Delta: \Gamma(X_{reg}, \Omega^2_{X_{reg}})(l) \to \Gamma (X_{reg}, \Omega^1_{X_{reg}})(l).$$ 
For $v \in \Gamma(X_{reg}, \Omega^2_{X_{reg}})(l)$, the Lie derivative $L_{\zeta}v$ of $v$ along $\zeta$
equals $l \cdot v$. If moreover $v$ is $d$-closed, then one has $l \cdot v = d(i_{\zeta}v)$ by 
the Cartan relation 
$$ L_{\zeta}v = d(i_{\zeta}v) + i_{\zeta}(dv). $$  
This means that $v$ is $d$-exact.  Q.E.D.  
\vspace{0.2cm}

{\bf  Proof of Theorem (3.1)}. Denote by $R$ the affine ring of $X$. 
By definition, $R$ has a natural grading 
$R = \oplus_{i \geq 0}R_i$ with $R_0 = \mathbf{C}$.
Let $j: X_{reg} \to X$ be the inclusion map. 
Since $j_*\Omega^2_{X_{reg}}$ is a coherent $O_X$-module, $M:= \Gamma (X_{reg}, \Omega^2_{X_{reg}})$ 
is a finitely generated, graded $R$-module: $M = \oplus M_i$. 
Each $M_i$ is a finite dimensional $\mathbf{C}$-vector space because 
$R_i = 0$ for $i<0$ and $R_0 = \mathbf{C}$. 
Our $\omega$ is an element of $M_l$ by the definition. 
Let $M_{l,{closed}}$ be the subspace of $M_l$ which consists of $d$-closed 
2-forms. Let $\mathrm{Aut}^{\mathbf{C}^*}(X)$ be the algebraic 
group of $\mathbf{C}^*$-equivariant automorphisms of $X$.
Then $\mathrm{Aut}^{\mathbf{C}^*}(X)$ acts on $M_{l,{closed}}$. 
Let $M_{l,closed}^0$ be the Zariski open subset of $M_{l, closed}$ which 
consists of non-degenerate 2-forms. In particular, $M_{l,closed}^0$ 
is connected. Since $\mathrm{Aut}^{\mathbf{C}^*}(X)$ preserves 
$M_{l,closed}^0$ as a set, we see that $M_{l,closed}^0$ is a single 
orbit of $\mathrm{Aut}^{\mathbf{C}^*}(X)$ by Proposition (3.2). Q.E.D.
\vspace{0.2cm}

{\bf Remark (3.3)}. Let $X$ be an affine variety defined by $f: = x^3 + y^3 + z^3 = 0$ in 
$\mathbf{C}^3$. Then $X$ has a natural $\mathbf{C}^*$-action with a fixed point 
$0 \in X$ and with $wt(x) = wt(y) = wt(z) = 1$. 
Then regular part $X_{reg}$ admits a symplectic 
form $\omega := Res(dx \wedge dy \wedge dz/f)$. 
The weight $l$ of $\omega$ is zero.
The blowing up of $X$ at $0$ gives 
us a resolution $\pi: \tilde{X} \to X$ with an exceptional curve $E$, which is an 
elliptic curve. The pull-back $\pi^*(\omega)$ is a meromorphic 2-form which has a 
pole along $E$. Thus $(X, \omega)$ is {\em not} a symplectic variety in the sense of [Be]. 
In this case, rigidity does not hold. In fact, $(X, t\cdot\omega)$ $(t \in \mathbf{C}^*)$ 
is a nontrivial Poisson deformation of $(X, \omega)$ (cf. [E-G]).  
We shall give here a short proof of this fact. By the argument of the proof of 
Rigidity Proposition, it suffices to prove that $\omega \in \Gamma (X_{reg}, \Omega^2_{X_{reg}})$ 
is not in the image of $d: \Gamma (X_{reg}, \Omega^1_{X_{reg}}) \to 
\Gamma (X_{reg}, \Omega^2_{X_{reg}})$.  Note that $\omega$ is a meromorphic 
2-form on $\tilde{X}$ having a pole along $E$ at order $1$. Thus one has 
$\omega \in \Gamma (\tilde{X}, \Omega^2_{\tilde{X}}(\mathrm{log} E)$. 
It can be checked that $\Gamma (\tilde{X}, \Omega^1_{\tilde X}(\mathrm{log}E)) 
\cong \Gamma (X_{reg}, \Omega^1_{X_{reg}})$. 
Let us consider the commutative diagram 

\begin{equation} 
\begin{CD} 
\Gamma (\tilde{X}, \Omega^1_{\tilde X}(\mathrm{log}E))  @>{Res}>>  \Gamma (E, O_E)\\ 
@V{d}VV @V{d}VV \\ 
\Gamma (\tilde{X}, \Omega^2_{\tilde{X}}(\mathrm{log} E)  @>{Res}>>  \Gamma (E, \Omega^1_E)   
\end{CD} 
\end{equation} 

Suppose that $\omega = d \eta$ for $\eta \in \Gamma (\tilde{X}, \Omega^1_{\tilde X}(\mathrm{log}E))$. 
Then one can write $$Res (\omega) = d Res (\eta)$$ by the commutative diagram. 
For any 1-cycle $\gamma$ on $E$, one has $$ \int_{\gamma} Res(\omega) = \int_{\gamma} d Res(\eta) = 0.$$ 
On the other hand, since $Res (\omega)$ is a nowhere vanishing 1-form on $E$, we 
should have $$\int_{\gamma} Res(\omega) \ne 0$$ for some 1-cycle 
$\gamma$ on $E$. This is a contradiction.          
\vspace{0.2cm}

{\bf Remark (3.4)}. Assume that $X$ has canonical singularities. 
Then the complex 
$$ \Gamma(X_{reg}, \Omega^1_{X_{reg}}) \stackrel{d}\to 
\Gamma(X_{reg}, \Omega^2_{X_{reg}}) \stackrel{d}\to 
\Gamma(X_{reg}, \Omega^3_{X_{reg}})$$ 
is exact. In particular, 
the complex 
$$ \Gamma(X_{reg}, \Omega^1_{X_{reg}})(0) \stackrel{d}\to 
\Gamma(X_{reg}, \Omega^2_{X_{reg}})(0) \stackrel{d}\to 
\Gamma(X_{reg}, \Omega^3_{X_{reg}})(0)$$  
is also exact. 
 
The proof goes as follows. Let $f: \tilde{X} \to X$ be a 
$\mathbf{C}^*$-equivariant resolution. Let $\alpha$ be a $d$-closed holomorphic  
2-form on $X_{reg}$.   
By [Na 2, Theorem 4] one has $\Gamma (\tilde{X}, \Omega^2_{\tilde X}) 
= \Gamma (X_{reg}, \Omega^2_{X_{reg}})$. Thus $f^*\alpha$ is a 
holomorphic 2-form on $\tilde{X}$. 
The following argument is based on [Na 3], \S 1. One can also find a 
similar argument in [Ka]. 

We first show that $[f^*\alpha] \in 
H^2(\tilde{X}, \mathbf{C})$ is zero. It is sufficient to prove that, for   
a small open neighborhood $U$ of $0 \in X$ (in the classical topology),  
$[f^*\alpha\vert_{f^{-1}(U)}] \in H^2(f^{-1}(U), \mathbf{C})$ 
is zero. In fact, the restriction map $H^2(\tilde{X}, \mathbf{C}) 
\to H^2(f^{-1}(U), \mathbf{C})$ is an isomorphism by the $\mathbf{C}^*$-action. 
One can blow up $\tilde{X}$ further to have a resolution $g: Z \to \tilde{X}$ so 
that $E:= (f \circ g)^{-1}(0)$ is a simple normal crossing divisor of $Z$. Since 
$R^1g_*\mathbf{C} = 0$, we have an injection $H^2(f^{-1}(U), \mathbf{C}) \to 
H^2((f \circ g)^{-1}(U), \mathbf{C})$. If we take $U$ sufficiently small, then 
$H^2((f \circ g)^{-1}(U), \mathbf{C}) \cong H^2(E, \mathbf{C})$. 
To prove that $[f^*\alpha\vert_{f^{-1}(U)}] = 0$ in 
$H^2(f^{-1}(U), \mathbf{C})$, we only have to check that 
$[(f \circ g)^*\alpha\vert_E] = 0$ in $H^2(E, \mathbf{C})$. 
Here we note that $H^2(E, \mathbf{C})$ has a mixed Hodge structure 
and $F^2(H^2(E, \mathbf{C}) = H^0(E, \hat{\Omega}^2_E)$. The sheaf $\hat{\Omega}^2_E$ is 
the quotient sheaf of $\Omega^2_E$ by the torsion subsheaf supported on 
$\mathrm{Sing}(E)$. Since $(f \circ g)^*\alpha\vert_E$ is a holomorphic 2-form 
on $E$, we have $[(f \circ g)^*\alpha\vert_E] \in H^0(E, \hat{\Omega}^2_E)$. 
But $H^0(E, \hat{\Omega}^2_E) = 0$ by [Na 3], Lemma (1.2). As a consequence, 
we have proved that $[f^*\alpha] \in H^2(\tilde{X}, \mathbf{C})$ is zero. 
 
Now look at the Hodge spectral sequence $$E^{p,q}_1 = H^q(\tilde{X}, \Omega^p_{\tilde X}) \Rightarrow 
H^{p+q}(\tilde{X}, \mathbf{C}).$$ 
Then $f^*\alpha = d\eta$ mod. $E^{0,1}_2$ with some holomorphic 1-form $\eta$ on $\tilde{X}$. 
Since $X$ has rational singularities, we have $E^{0,1}_1 = 0$; hence $E^{0,1}_2 = 0$.  
Thus $f^*\alpha = d\eta$. This clearly shows that the original complex is exact.   
\vspace{0.2cm}

{\bf Proposition (3.5)} (Strong rigidity). {\em 
Assume, in addition, that $X$ has canonical singularities. 
Let $(X_1, \{\;,\;\}_{\epsilon})$ be a $\mathbf{C}^*$-equivariant Poisson deformation of 
$(X, \{\;,\;\})$ over $S_1$ in such a way that $\mathbf{C}^*$ acts on $S_1$ trivially. 
Then this Poisson deformation is a trivial one.} 
\vspace{0.2cm}

{\em Proof}. The difference from Proposition (3.2) is that we do not assume 
that $X_1 = X \times S_1$.   
Let $f: \mathcal{X} \to \mathbf{A}^d$ be a $\mathbf{C}^*$-equivariant 
universal Poisson deformation of 
$X$ over an affine space $\mathbf{A}^d$ constructed in [Na 1]. 
Note that there is a Poisson isomorphism $\iota: X \cong f^{-1}(0)$.    
The $\mathbf{C}^*$-action on $X$ induces a $\mathbf{C}^*$-action on 
the base space $\mathbf{A}^d$ of the universal Poisson deformation. 
By the construction of $f$ ([ibid]), this action has only positive weights. 

The infinitesimal Poisson deformation $X_1 \to S_1$ determines a map 
$S_1 \to \mathbf{A}^d$ which sends the closed point of $S_1$ to the origin of 
$\mathbf{A}^d$. Assume that this is a closed immersion; namely, $S_1 \subset 
\mathbf{A}^d$. By the assumption, the $\mathbf{C}^*$-action on $\mathbf{A}^d$ 
restricts to the trivial action on $S_1$. This contradicts that the 
$\mathbf{C}^*$-action on $\mathbf{A}^d$ has only positive weights. 
Thus the map $S_1 \to \mathbf{A}^d$ is a constant map.  Q.E.D.   
\vspace{0.2cm}
       
{\bf 4. Projectivised cone and contact structures} 

In this section $(X, \omega)$ is a pair of a normal affine variety $X$ of dimension $2d$ with a good $\mathbf{C}^*$-action and an algebraic symplectic 2-form $\omega$ on $X_{reg}$ with 
{\em positive} weight $l$.

(4.1) {\em Projectivised cone} 
  
Let $R$ be the affine ring of $X$. By definition, $R$ has a natural grading 
$R = \oplus_{i \geq 0}R_i$ with $R_0 = \mathbf{C}$. We put $$\mathbf{P}(X) := \mathrm{Proj}(\oplus_{i \geq 0}R_i).$$    
Let $x_0$, $x_1$, ..., $x_n$ be homogeneous minimal generators of the $\mathbf{C}$-algebra 
$R$ and put $a_i = wt(x_i)$. 
Then $\mathbf{P}(X)$ is naturally embedded in the weighted projective space 
$\mathbf{P}(a_0, a_1, ..., a_n)$. Let $V \to \mathbf{C}^{n+1}$ be the weighted blowing up of 
$\mathbf{C}^{n+1}$ with weight $(a_0, ..., a_n)$. Then the fibre over the origin $0 \in \mathbf{C}^{n+1}$ 
is isomorphic to $\mathbf{P}(a_0, ..., a_n)$. 
The singular locus $\mathrm{Sing}$ of $V$ is contained in the fibre over the origin; hence 
one can regard $\mathrm{Sing}(V)$ as a subset of $\mathbf{P}(a_0, ..., a_n)$. 
In this identification $\mathrm{Sing}(V)$ is the locus where the projection map 
$\mathbf{C}^{n+1}-\{0\} \to \mathbf{P}(a_0, ..., a_n)$ is not a $\mathbf{C}^*$-bundle. 
As the subsets of $\mathbf{P}(a_0, ..., a_n)$, we may take the intersection of 
$\mathbf{P}(X)$ and $\mathrm{Sing}(V)$.  We assume that 
$\mathbf{P}(X) \cap \mathrm{Sing}(V)$ has codimension at least $2$ in $\mathbf{P}(X)$. 
Notice that this condition is not necessarily satisfied. For example, 
$A_m$-surface singularity $x_0^2 + x_1^2 + x_2^{m+1} = 0$ for $m \geq 2$ 
does not satisfy this condition.


Let $\mathbf{P}(X)^0$ be the open subset obtained by 
excluding this subset and $\mathrm{Sing}(\mathbf{P}(X))$ from $\mathbf{P}(X)$. 
Note that $\mathrm{Codim}_{\mathbf{P}(X)}(\mathbf{P}(X) - 
\mathbf{P}(X)^0) \geq 2$.    
There is a natural projection 
$$p: X - \{0\} \to \mathbf{P}(X),$$ which is a $\mathbf{C}^*$-fibration and is actually a  
$\mathbf{C}^*$-fibre bundle over $\mathbf{P}(X)^0$. 
We put $X^0 := p^{-1}(\mathbf{P}(X)^0)$. Let $O(1)$ be the tautological sheaf on $\mathbf{P}(a_0, ..., a_n)$ 
and put $O_{\mathbf{P}(X)}(1) := O(1) \otimes_{O_{\mathbf{P}(a_0, ..., a_n)}}O_{\mathbf{P}(X)}$. 
Then $O_{\mathbf{P}(X)}(1)\vert_{\mathbf{P}(X)^0}$ is an invertible sheaf on $\mathbf{P}(X)^0$.    
Let $L \in \mathrm{Pic}(\mathbf{P}(X)^0)$ be the corresponding line bundle to this sheaf. More exactly, 
$O_{\mathbf{P}(X)}(1)\vert_{\mathbf{P}(X)^0}$ is the sheaf of sections of $L$. Denote by $L^{-1}$ the dual 
line bundle of $L$ and denote by $(L^{-1})^{\times}$ the  
$\mathbf{C}^*$-bundle which is obtained from $L^{-1}$ by removing the zero section.  Then $X^0$ coincides 
with $(L^{-1})^{\times}$ and the natural projection $$\pi: (L^{-1})^{\times} \to \mathbf{P}(X)^0$$ coincides 
with $p\vert_{X^0}$.  Note that there is a canonical trivialisation $$\pi^*L \cong O_{(L^{-1})^{\times}}.$$ 
Recall that $l$ is the weight of $\omega$.
Later we will use the trivialisation $$ \pi^*(L^{\otimes l}) \to O_{(L^{-1})^{\times}}$$ induced by this canonical trivialisation.   
\vspace{0.2cm}

(4.2) {\em Contact structure on a complex manifold} 

We shall briefly review a contact complex manifold according to LeBrun [LeB]. 
Let $Z$ be a complex manifold of dimension $2d+1$. A contact structure on $Z$ is an exact sequence of vector bundles 
$$ 0 \to D \to TZ \stackrel{\theta}\to M \to 0,$$ with $\mathrm{rank}(D) = 2d$ and $\mathrm{rank}(M) = 1$ so 
that $d\theta\vert_D$ induces a non-degenerate pairing on $D$. 
By using the formula for exterior derivation 
$$ d\theta (x,y) = x(\theta (y)) - y(\theta (x))  - \theta ([x,y]) $$ one can check that this is equivalent 
to saying that $[\;, \;] : D \times D \to TZ/D (=M)$ is non-degenerate. 
We call $M$ the contact line bundle.   
As is well known, infinitesimal automorphisms of $Z$ are controlled by the cohomology group $H^0(Z, \Theta_Z)$. 
An infinitesimal automorphism of $Z$ is said to be contact if it preserves the contact structure. 

{\bf Proposition (4.2.1)} ([LeB, Proposition 2.1]) {\em Let
$$ 0 \to O(D) \to \Theta_Z \stackrel{\theta}\to O(M) \to 0 $$ be 
the exact sequence of sheaves determined by the contact structure. Then there 
is a map $s: O(M) \to \Theta_Z$ of $\mathbf{C}$-modules (not of $O_Z$-modules) that splits the sequence 
above and the group of infinitesimal contact automorphisms coincides with 
$s(H^0(Z, O(M))$. }    
\vspace{0.2cm}

{\bf Corollary (4.2.2)} (cf. [ibid, Proposition 2.2]) {\em Fix a line bundle $M$ on $Z$. Assume that $TZ \stackrel{\theta}\to M$ is a 
contact structure on $Z$ so that $H^0(Z, O(D)) = 0$. Then $\theta$ is a unique contact structure 
with contact line bundle $M$. }
\vspace{0.2cm}

(4.3) {\em Quasi-contact structure on $\mathbf{P}(X)$} 

One can generalise the notion of 
contact structure to a singular variety. Let $Z$ be a normal variety. 
Here a quasi-contact structure\footnote{We do not assume that $j_*M$ is a 
line bundle on $Z$. As we will define in (4.4), if $j_*M$ is a line bundle on $Z$, we call 
it a contact structure on $Z$.} on $Z$ is just a contact structure on an open set $Z^0 \subset Z_{reg}$ 
with $\mathrm{codim}_Z (Z - Z^0) \geq 2$. By the definition,             
there are a line bundle $M$ on $Z^0$ and a vector bundle $D$ on $Z^0$ of rank $2d$ which 
fit into an exact sequence $$ 0 \to O(D) \to \Theta_{Z^0} \to O(M) \to 0.$$ 
Since the degeneracy locus of a contact form has codimension one, a contact structure 
on $Z^0$ uniquely extends to that on $Z_{reg}$. Thus we may say that a qusi-contact structure on 
$Z$ is a contact structure on $Z_{reg}$.  
Let $j: Z^0 \to Z$ be the natural inclusion map. Then we have an exact sequence 
$$ 0 \to j_*O(D) \to \Theta_Z \to j_*O(M) \to 0. $$ 
Note that the last map is surjective by Proposition (4.2.1).  
\vspace{0.2cm}

Let us return to the original situation.  The complement of $\mathbf{P}(X)^0$ 
in $\mathbf{P}(X)$ has at least codimension $2$.         
Let us introduce a quasi-contact structure on $\mathbf{P}(X)$. This is a slight modification 
of the argument in [LeB, p.425], where the case $l = 1$ is treated. 
Recall that we have a $\mathbf{C}^*$-bundle $p\vert_{X^0}: X^0 \to \mathbf{P}(X)^0$ 
and it is identified with $\pi: (L^{-1})^{\times} \to \mathbf{P}(X)^0$. 

For $\theta \in H^0(\mathbf{P}(X)^0, \Omega^1_{\mathbf{P}(X)^0}(L^{\otimes l}))$, 
the pull-back $\pi^*(\theta)$ is regarded as an element of 
$H^0((L^{-1})^{\times}, \Omega^1_{(L^{-1})^{\times}})$ by the trivialisation  
$\pi^*(L^{\otimes l}) \to O_{(L^{-1})^{\times}}$. 
 
By the assumption we have a symplectic 2-form $\omega$ on $(L^{-1})^{\times}$ with weight $l$. 
As a $\mathbf{C}^*$-bundle, there is a natural $\mathbf{C}^*$-action on $(L^{-1})^{\times}$.
Let $\zeta$ be the vector field which generates the $\mathbf{C}^*$-action.  
Then one can write  $\omega(\zeta, \cdot) = \pi^*\theta$ with an element $\theta 
\in H^0(\mathbf{P}(X)^0, \Omega^1_{\mathbf{P}(X)^0}(L^{\otimes l}))$. This 
$\theta$ gives a contact structure on $\mathbf{P}(X)^0$ with contact line bundle 
$L^{\otimes l}$. Conversely, if a contact structure $\theta \in 
H^0(\mathbf{P}(X)^0, \Omega^1_{\mathbf{P}(X)^0}(L^{\otimes l}))$ is given 
to $\mathbf{P}(X)^0$, then $d\pi^*(\theta)$ becomes a holomorphic symplectic 
2-form on $(L^{-1})^{\times}$ with weight $l$. Note that we need the assumption $l \ne 0$ 
to get the correspondence between symplectic structures of weight $l$ and contact structures.   
\vspace{0.2cm}
   
(4.4) {\em Contact orbifold structure and Jacobi orbifold structure} 

In (4.1) we imposed a rather technical assumption; namely 
$\mathbf{P}(X) \cap \mathrm{Sing}(V)$ has at least 
codimension $2$ in $\mathbf{P}(X)$. 

In the remainder of this section we do {\em not} assume this. 

In a general case 
a possible structure would be a {\em contact orbifold structure}. 
Let us consider a normal variety $Z$ and 
a line bundle $M$ on $Z$.  
A contact structure on $Z$ (with contact line bundle $M$) is a contact structure on the  
Zariski open set $Z_{reg}$ (as a complex manifold) with contact line bundle 
$M\vert_{Z_{reg}}$. A contact form $\theta$ is regarded as a section of 
$\underline{\mathrm{Hom}}(\Theta_Z, M)$.   
A contact orbifold $Y$ is a normal variety with 
the following data: $Y = \cup U_{\alpha}$ is an open covering of $Y$ and, for each 
$\alpha$, there is a finite Galois covering $\varphi_{\alpha}: \tilde{U}_{\alpha} \to 
U_{\alpha}$ such that the (possibly singular but normal) variety $\tilde{U}_{\alpha}$ admits a line bundle 
$M_{\alpha}$ and a contact form $\theta_{\alpha}$ with contact line bundle $M_{\alpha}$.  
These data should satisfy a compatibility condition.  
If $U_{\alpha} 
\cap U_{\beta} \ne \emptyset$, then we form a diagram 
$$ \tilde{U}_{\alpha} \stackrel{p_{\alpha}}\leftarrow \tilde{U}_{\alpha} \times_Y 
\tilde{U}_{\beta} \stackrel{p_{\beta}} \to \tilde{U}_{\beta}.$$ 
Let $(\tilde{U}_{\alpha} \times_Y \tilde{U}_{\beta})^n$ be the normalisation 
of $\tilde{U}_{\alpha} \times_Y \tilde{U}_{\beta}$.  
Denote by $p_{\alpha}^n$ the composite of the normalisation map and $p_{\alpha}$. 
We then assume that $p_{\alpha}^n$ and $p_{\beta}^n$ are both etale. 
Moreover, as the compatibility condition we assume that there is an isomorphism of line bundles 
$$g_{\beta, \alpha}: (p_{\alpha}^n)^*M_{\alpha} \to (p_{\beta}^n)^*M_{\beta}$$ and that 
$$(p_{\alpha}^n)^*(\theta_{\alpha}) =  (p_{\beta}^n)^*(\theta_{\beta}).$$ 
Finally, for any $\alpha$, $\beta$ and $\gamma$ with $U_{\alpha} \cap U_{\beta} 
\cap U_{\gamma} \ne \emptyset$, we should have 
$$ g_{\alpha, \beta} \circ g_{\beta, \gamma} \circ g_{\gamma, \alpha} = id $$ 
on $$(\tilde{U}_{\alpha} \times_Y \tilde{U}_{\beta} \times_Y \tilde{U}_{\gamma})^n.$$
In other words, $\{M_{\alpha}\}$ is an orbifold line bundle $\mathcal{M}$ on $Y^{orb}$, and 
$\{\theta_{\alpha}\}$ is a global section of $\underline{\mathrm{Hom}}(\Theta_{Y^{orb}}, 
\mathcal{M})$.

The most natural structure would be actually a {\em Jacobi structure} ([Li]).  
This is very similar to the fact that a Poisson structure would be more natural than a 
symplectic structure in the singular case.   
If a normal variety has a contact structure in the sense above, then we have a pairing map 
$$ O(M)\vert_{Z_{reg}} \times O(M)\vert_{Z_{reg}} \to O(M)\vert_{Z_{reg}} $$ 
defined by $(u,v) \to \theta ([s(u), s(v)])$. Here $s$ is the map defined in Proposition in {\bf 6}. 
By the normality this pairing uniquely extends to 
$$\{\;, \;\}: O(M) \times O(M) \to O(M).$$ The bracket satisfies the Jacobi identity, but it is 
no more a bi-derivation. We call it a Jacobi structure on $Z$.  
The Jacobi structure is generalised to orbifold version in a similar way as the contact orbifold structure 
was defined. A contact orbifold structure determines a Jacobi orbifold structure.  
\vspace{0.2cm}

{\bf Theorem (4.4.1)}. {\em The projectivised cone   $\mathbf{P}(X)$ has a 
contact orbifold structure.} 
\vspace{0.2cm}

{\em Proof}. First note that $\mathbf{P}(a_0, ..., a_n)$ has a natural orbifold structure. 
In fact, let $\mathbf{C}^{n+1}-\{0\} \to \mathbf{P}(a_0, ..., a_n)$ be the quotient map of the 
$\mathbf{C}^*$-action $(x_0, ..., x_n) \to (t^{a_0}x_0, ..., t^{a_n}x_n)$.  
Restrict this map to $W_i := \{x_i = 1\} \subset \mathbf{C}^{n+1}$. Then one has a map 
$W_i \to \mathbf{P}(a_0, ..., a_n)$ for each $i$ and these maps give an orbifold structure of 
$\mathbf{P}(a_0, ..., a_n)$. We show that $\mathbf{P}(a_0, ..., a_n)$ admits an orbifold 
line bundle $O_{\mathbf{P}(a_0, ..., a_n)}(1)$. There is a finite Galois cover 
$$\mathbf{P}(a_0, ..., a_{i-1}, 1, a_{i+1}, ..., a_n) \to \mathbf{P}(a_0, ..., a_n)$$ 
defined by $$(x_0, ..., x_n) \to (x_0, ..., x_i^{a_i}, ..., x_n)$$ 
for each $i$. One can identify $W_i$ with the open set of 
$\mathbf{P}(a_0, ..., a_{i-1}, 1, a_{i+1}, ..., a_n)$ defined by $x_i \ne 0$. 
Let $$\tilde{L}_i := O_{\mathbf{P}(a_0, ..., a_{i-1}, 1, a_{i+1}, ..., a_n)}(1)\vert_{W_i}.$$
Then $\{\tilde{L}_i\}_{0 \le i \le n}$ give an orbifold line bundle on $\mathbf{P}(a_0, ..., a_n)$. 
In fact, the $\mathbf{Z}/a_0\mathbf{Z} \times ... \times \mathbf{Z}/a_n\mathbf{Z}$-Galois cover 
$$ \mathbf{P}(1, ..., 1) \to \mathbf{P}(a_0, ..., a_n)$$ 
is a smooth global cover (cf. [Mu], Section 2) in the sense that it is factorized as 
$$\mathbf{P}(1, ..., 1) \to \mathbf{P}(a_0, ..., a_{i-1}, 1, a_{i+1}, ..., a_n) \to 
\mathbf{P}(a_0, ..., a_n)$$ for any $i$. The tautological line bundle 
$O_{\mathbf{P}(1, ..., 1)}(1)$ has a  
$G := \mathbf{Z}/a_0\mathbf{Z} \times ... \times \mathbf{Z}/a_n\mathbf{Z}$ linearization 
defined by $x_i \to \zeta_i^{m_i}x_i$ $(0 \le i \le n)$ for a primitive $a_i$-th root $\zeta_i$ of unity 
and $m_i \in \mathbf{Z}/a_i\mathbf{Z}$. 
Then $O_{\mathbf{P}(1, ..., 1)}(1)\vert_{x_i \ne 0}$ with the action of $G_i := 
{\mathbf{Z}/a_i\mathbf{Z}}$ is the pullback of $\tilde{L}_i$. 
This is equivalent to giving an orbifold line bundle of $\mathbf{P}(a_0, ..., a_n)$ 
(cf. [Ibid]). 
The merit of introducing the orbifold structure is the following.
Let $\Sigma \subset \mathbf{P}(a_0, ..., a_n)$ be the union\footnote{If we identify 
$\mathbf{P}(a_0, ..., a_n)$ with the central fibre of the weighted blowing up 
$V \to \mathbf{C}^{n+1}$, then $\Sigma$ coincides with 
$\mathrm{Sing}(V)$.} of the ramification loci of 
the coverings $W_i \to \mathbf{P}(a_0, ..., a_n)$.  
Each fibre of the projection map $\mathbf{C}^{n+1} - \{0\} \to \mathbf{P}(a_0, ..., a_n)$ 
is isomorphic to $\mathbf{C}^*$, but the fibres over the points contained in $\Sigma$ are 
multiple fibres. However, if we take the normalisation 
$(W_i \times_{\mathbf{P}(a_0, ..., a_n)}(\mathbf{C}^{n+1}-\{0\}))^n$ 
of the fibre product of $W_i$ and $\mathbf{C}^{n+1}-\{0\}$ over $\mathbf{P}(a_0, ..., a_n)$, then 
the first projection 
$$(W_i \times_{\mathbf{P}(a_0, ..., a_n)}(\mathbf{C}^{n+1}-\{0\}))^n \to W_i$$ is a 
$\mathbf{C}^*$-bundle and the second projection 
$$(W_i \times_{\mathbf{P}(a_0, ..., a_n)}(\mathbf{C}^{n+1}-\{0\}))^n \to \mathbf{C}^{n+1}-\{0\}$$ 
is an etale map.

Put $U_i := W_i \cap X$ and $L_i := \tilde{L}_i\vert_{U_i}$.  
Then an orbifold structure of $\mathbf{P}(X)$ is given by $\{U_i \to \mathbf{P}(X)\}$. 
Moreover $\{L_i\}$ give an orbifold line bundle $\mathcal{L}$ on $\mathbf{P}(X)$. 
The orbifold line bundle $\mathcal{L}$ is called the tautological line bundle. 
Let $\mathcal{M}$ be the orbifold line bundle on $\mathbf{P}(X)$ defined by 
$\mathcal{L}^{\otimes l}$. 

Let $X_i$ be the normalisation of the fibre product $U_i \times_{\mathbf{P}(X)}(X - \{0\})$. 
Then the first projection $X_i \to U_i$ is a $\mathbf{C}^*$-bundle and the second 
projection $X_i \to X - \{0\}$ is an etale map. 
Let $\omega_i$ be the pullback of $\omega$ by the map $(X_i)_{reg} \to X_{reg}$. 
As in (4.3), $\omega_i$ defines a contact structure on $(U_i)_{reg}$ with contact line 
bundle $L_i^{\otimes l}\vert_{(U_i)_{reg}}$. These contact structures are glued together to give 
a contact orbifold structure on $\mathbf{P}(X)$ with contact line bundle 
$\mathcal{M}$. Q.E.D. 
\vspace{0.2cm}
 
We shall briefly recall the cohomology for an orbifold. Let $\pi_i: U_i \to \mathbf{P}(X)$ 
be the orbifold charts with $G_i = \mathrm{Gal}(\pi_i)$ and let 
$p^n_{i,j}: (U_i \times_{\mathbf{P}(X)}U_j)^n \to U_i$ be the 
projection maps from the normalisation of $U_i \times_{\mathbf{P}(X)}U_j$ to $U_i$. 
An orbifold $O_{\mathbf{P}(X)}$-module 
$F$ is a collection $\{F_i\}$ of $O_{U_i}$-modules glued together by 
$$g_{i,j}: (p^n_{i,j})^*F_i \cong (p^n_{j,i})^*F_j$$  
compatible on the triple overlaps. 
If we put, in particular, $j = i$, then this means that $F_i$ has a $G_i$-linearisation. 
The sheaves $(\pi_i)^{G_i}_*{F_i}$ of $G_i$-invariant sections of $F_i$ are glued together to 
give a sheaf $\bar{F}$ on $\mathbf{P}(X)$. The space $\Gamma (\mathbf{P}(X), F)$ of global 
sections of $F$ is nothing but $\Gamma (\mathbf{P}(X), \bar{F})$. We define $H^p(\mathbf{P}(X), F) 
:= H^p(\mathbf{P}(X), \bar{F})$. 
When $F_i$ are all invertible sheaves, $F$ is called an orbifold line bundle. 
Even if $F$ is an orbifold line bundle, $\bar{F}$ is not necessarily a line bundle. 
Let $\mathrm{Pic}(\mathbf{P}(X)^{orb})$ be the group of isomorphism classes of orbifold 
line bundles on $\mathbf{P}(X)$. In general $\mathrm{Pic}(\mathbf{P}(X)^{orb})$ is not 
isomorphic to $H^1(\mathbf{P}(X), O^*_{\mathbf{P}(X)})$. In order to capture the orbifold 
line bundles we consider the \v{C}ech complex 
$$ \prod_i \Gamma (U_i, O^*_{U_i}) \to \prod_{i,j} \Gamma (U_{i,j}, O^*_{U_{i,j}}) \to \prod_{i,j,k} 
\Gamma (U_{i,j,k}, O^*_{U_{i,j,k}}) \to ... $$
We denote by $H^p_{orb}(\mathcal{U}, O^*_{\mathbf{P}(X)})$  
the $p$-th cohomology of this complex. The inductive limit $\lim H^p_{orb}(\mathcal{U}, O^*_{\mathbf{P}(X)})$ 
for the admissible orbifold charts is the $p$-th \v{C}ech orbifold cohomology 
$H^p_{orb}(\mathbf{P}(X), O^*_{\mathbf{P}(X)})$ of $O^*_{\mathbf{P}(X)}$.  
Then $$\mathrm{Pic}(\mathbf{P}(X)^{orb}) \cong H^1_{orb}(\mathbf{P}(X), O^*_{\mathbf{P}(X)}).$$ 
As in the proof of (4.4.1), let $\mathcal{L}$ be the tautological line bundle on $\mathbf{P}(X)$.   

We put $$ R := \oplus_{n \geq 0} H^0(\mathbf{P}(X), \mathcal{L}^{\otimes n}).$$ 
Then $X = \mathrm{Spec}(R)$. If we pull back the projection map $p: X - \{0\} \to \mathbf{P}(X)$ 
by $U_i \to \mathbf{P}(X)$ and take the normalisation, then we have a $\mathbf{C}^*$-bundle 
$X_i \to U_i$ and the induced map $X_i \to X -\{0\}$ is etale. The contact structure $\theta_i$ 
induces a symplectic structure $\omega_i$ of weight $l$ on $X_i$. These symplectic structures 
$\{\omega_i\}$ descend to a symplectic structure $\omega$ of weight $l$ on $X$.  
In this way, $(X, \omega)$ is recovered from the 
contact structure $(\mathcal{M}, \theta)$. 
\vspace{0.2cm}

{\bf Example (4.4.2)}. 

(i) Let us consider a Du Val singularity $X$ (of type $A_n$, $D_n$ ($n \geq 4$), or $E_n$ ($n = 6,7,8$) 
with a symplectic structure $\omega$ of weight 2. By Theorem (4.4.1) $\mathbf{P}(X)$ has a contact orbifold 
structure. Du Val singularities of different type determine mutually different contact orbifold 
structures. But the underlying variety $\mathbf{P}(X)$ are all $\mathbf{P}^1$. 
In other words, $\mathbf{P}^1$ has infinitely many different contact orbifold structures. 

(ii) The odd dimensional projective space $\mathbf{P}^{2n+1}$ has a contact structure with a 
contact line bundle $M = O(2)$. We have two choices of the tautological line bundle 
$L$: $L =O(1)$ or $L = O(2)$. The weights $l$ are then respectively $2$ and $1$. 
The corresponding symplectic variety $(X, \omega)$ is isomorphic to $\mathbf{C}^{2n+2}$ with 
the standard symplectic structure in the first case.     
In the second case $(X, \omega)$ is isomorphic to $\mathbf{C}^{2n+2}/\{+1, -1\}$ with  
a symplectic 2-form $\omega_0$. Here $-1$ acts on $\mathbf{C}^{2n+2}$ by $x_i \to -x_i$ 
$(1 \le i \le 2n+2)$ and $\omega_0$ is the symplectic structure induced from the 
standard symplectic form $dx_1 \wedge dx_2 + ... + dx_{2n+1} \wedge dx_{2n+2}$.  
\vspace{0.2cm}

{\bf 5. Rigidity of contact orbifold structures} 

(5.1)  Let $(X, \omega)$ be a pair of a normal affine variety $X$  of dimension $2d$ with a good $\mathbf{C}^*$-action and an algebraic symplectic 2-form $\omega$ on $X_{reg}$ with 
positive weight $l$. In the previous section we have attached a contact orbifold structure $(\mathcal{M}, \theta)$ to the projectivised cone $\mathbf{P}(X)$. In this section we consider 
the deformation of such a contact orbifold. 
  
Let $(Y, \{U_{\alpha}\}, \{M_{\alpha}\}, \{\theta_{\alpha}\})$ be a contact orbifold and let $S$ be a punctured space with $0 \in S$. 
A flat deformation of $(Y, \{U_{\alpha}\}, \{M_{\alpha}\}, \{\theta_{\alpha}\})$ is a 
a flat surjective map $\mathcal{Y} \to S$  
together with the covering charts $\Phi_{\alpha}: \tilde{\mathcal U}_{\alpha} \to \mathcal{Y}$ 
such that 

(i) $\mathcal{Y}$ is a flat deformation of $Y$, 

(ii) for each $\alpha$, $\Phi_{\alpha}$ is a Galois covering with the Galois 
group $G_{\alpha}$ which is a flat deformation of $\varphi_{\alpha}: \tilde{U}_{\alpha} \to Y$ over $S$ 
so that the maps 
$$(\tilde{\mathcal U}_{\alpha} \times_{\mathcal Y} \tilde{\mathcal U}_{\beta})^n \to 
\tilde{\mathcal U}_{\alpha}$$ are etale.    

(iii) there are line bundles $\mathcal{M}_{\alpha}$ on $\tilde{\mathcal U}_{\alpha}$ 
that restrict to give $M_{\alpha}$ on $\tilde{U}_{\alpha}$, and $\{\mathcal{M}_{\alpha}\}$ 
are glued together to give an orbifold line bundle of $\mathcal{Y}$, and finally     

(iv) all contact structures $\theta_{\alpha}$ on $\tilde{U}_{\alpha}$ extend compatibly to  
the contact structures $\Theta_{\alpha}$ on $\tilde{\mathcal U}_{\alpha}$ with the contact 
line bundles $\mathcal{M}_{\alpha}$.    
 
Of course, one can start from a different admissible orbifold charts $\{U'_{\alpha}\}$ of $Y$ and consider its 
flat deformation. A flat deformation of the contact orbifold structure $(Y, \mathcal{M}, \theta)$ is exactly 
an equivalence classe of that of $(Y, \{U_i\}, \{M_{\alpha}\}, \{\theta_{\alpha}\})$.  
\vspace{0.2cm}

Now let us return to the contact orbifold $(\mathbf{P}(X), \mathcal{M}, \theta)$. 
\vspace{0.2cm}
 
{\bf Proposition (5.2)}. {\em Assume that $X$ has only canonical singularities. 
Then the contact orbifold structure $(\mathbf{P}(X), \mathcal{M}, \theta)$ is 
rigid under a small flat deformation.}
\vspace{0.15cm}

This is a counterpart of Proposition (3.5) in the contact geometry. 
\vspace{0.15cm}

{\em Proof}. Recall the construction of the contact orbifold structure (cf. Proof of 
Theorem (4.4.1)). With the same notation as in the proof of (4.4.1), the map 
$X_i \to X -\{0\}$ is an etale map. Since $X$ has canonical singularities, $X_i$ also 
has canonical singularities. Since $X_i$ is a $\mathbf{C}^*$-bundle over $U_i$, 
we see that $U_i$ has canonical singularities; hence all orbifold charts of 
$\mathbf{P}(X)$ have canonical singularities.   

Let $(Y, \mathcal{M}_1, \theta_1)$ be an infinitesimal deformation of   
$(\mathbf{P}(X), \mathcal{M}, \theta)$ over $S_1 := \mathrm{Spec} \mathbf{C}[\epsilon]$. 
The orbifold line bundle $\mathcal{L}$ defines an element $[\mathcal{L}]$ of 
$H^1_{orb}(\mathbf{P}(X), O^*_{\mathbf{P}(X)})$. 
Note that there is an exact sequence  
$$ 0 \to O_{\mathbf{P}(X)} \to O^*_Y \to O^*_{\mathbf{P}(X)} \to 1,$$ 
where the map $O_{\mathbf{P}(X)} \to O^*_Y$ is defined by $f \to 1 + \epsilon\cdot f$. 
This exact sequence yields the exact sequence 
$$ H^1_{orb}(\mathbf{P}(X), O_{\mathbf{P}(X)}) \to H^1_{orb}(Y, O^*_Y) \to 
H^1_{orb}(\mathbf{P}(X), O^*_{\mathbf{P}(X)}) \stackrel{\delta}\to H^2_{orb}(\mathbf{P}(X), O_{\mathbf{P}(X)}).$$ 
Note that $\mathcal{M} = \mathcal{L}^{\otimes l}$. Since $\mathcal{M}$ extends to $\mathcal{M}_1$, 
one has $\delta([\mathcal{M}]) = 0$. This means that $l\cdot \delta([\mathcal{L}]) = 0$. 
As $H^2_{orb}(\mathbf{P}(X), O_{\mathbf{P}(X)})$ is a $\mathbf{C}$-vector space, 
one has $\delta([\mathcal{L}]) = 0$. Thus one can find an orbifold line bundle $\mathcal{L}_1$ on $Y$ that  
is an extension of $\mathcal{L}$. Then $l\cdot [\mathcal{L}_1] - [\mathcal{M}_1] \in H^1_{orb}(Y, O^*_Y)$ 
is the image of an element $\eta \in H^1_{orb}(\mathbf{P}(X), O_{\mathbf{P}(X)})$ by the map 
$H^1_{orb}(\mathbf{P}(X), O_{\mathbf{P}(X)}) \to H^1_{orb}(Y, O^*_Y)$. If we replace $\mathcal{L}_1$ by the 
orbifold line bundle corresponding to $[\mathcal{L}_1] - 1/l\cdot \eta$, then we may assume that 
$\mathcal{L}_1^{\otimes l} = \mathcal{M}_1$.    
The exact sequences  
$$ 0 \to \epsilon\cdot \mathcal{L}^{\otimes n} \to \mathcal{L}_1^{\otimes n} 
\to \mathcal{L}^{\otimes n} \to 0$$ 
yield the exact sequences
$$ H^0(Y, \mathcal{L}_1^{\otimes n}) \to  H^0(\mathbf{P}(X), \mathcal{L}^{\otimes n}) 
\to H^1(\mathbf{P}(X), \mathcal{L}^{\otimes n}).$$ 
By the next lemma, we see that the maps 
$H^0(Y, \mathcal{L}_1^{\otimes n}) \to  H^0(\mathbf{P}(X), \mathcal{L}^{\otimes n})$ 
are all surjective for $n \geq 0$. 
We put 
$$ \mathcal{R} := \oplus_{n \geq 0} H^0(Y, \mathcal{L}_1^{\otimes n}).$$
and define $\mathcal{X} := \mathrm{Spec}(\mathcal{R})$. 
Then $\mathcal{X}$ is an infinitesimal deformation of $X$ over $S_1$. 
Moreover the contact structure $(\mathcal{M}_1, \theta_1)$ of $Y$ defines a symplectic 
structure $\omega_1$ on $\mathcal{X}$. As a consequence, we have obtained an infinitesimal 
deformation $(\mathcal{X}, \omega_1)$ of $(X, \omega)$.   
By the construction $(\mathcal{X}, \omega_1)$ has a $\mathbf{C}^*$-action. If one regards 
$S_1$ as a $\mathbf{C}^*$-space with trivial action, then the map  
$\mathcal{X} \to S_1$ is $\mathbf{C}^*$-equivariant. By Proposition (3.5)  
$(\mathcal{X}, \omega_1)$ is a trivial deformation of $(X, \omega)$. 
In particular, $(Y, \mathcal{M}_1, \theta_1)$ is also a trivial deformation of    
$(\mathbf{P}(X), \mathcal{M}, \theta)$.  
\vspace{0.2cm}   

{\bf Lemma (5.2.1)}. {\em $$H^1(\mathbf{P}(X), \mathcal{L}^{\otimes n}) = 0$$ 
for all $n \geq 0$.}  
\vspace{0.2cm}

{\em Proof}. As remarked above, each orbifold of $\mathbf{P}(X)$ has rational 
Gorenstein singularities. 
Since $\mathbf{P}(X)$ is locally the quotient variety of an orbifold 
chart by a finite group, the log variety $(\mathbf{P}(X), 0)$ (with the zero boundary divisor) 
has log terminal singularities by 
[Kaw, Proposition 1.7]\footnote{The map $\pi: V \to X$ in Proposition 1.7 is 
assumed to be etale in codimension 1. But this condition is not necessary to prove the 
``if part".} .  
Moreover $\mathbf{P}(X)$ has a contact orbifold structure; thus 
$K^{orb}_{\mathbf{P}(X)} = \mathcal{M}^{-d-1}$ if $\dim \mathbf{P}(X) = 2d + 1$. 
Since $\mathcal{M} = \mathcal{L}^{\otimes l}$ is ample, $K^{orb}_{\mathbf{P}(X)}$ 
is negative.

Let $\pi_i : U_i \to \mathbf{P}(X)$ be an orbifold chart. Then $U_i \to \pi_i(U_i)$ is a   
$\mathbf{Z}/a_i\mathbf{Z}$- Galois cover.  
One can write $K_{U_i} = \pi_i^*K_{\mathbf{P}(X)} + B_i$ with an effective divisor 
$B_i$ on $U_i$ whose support coincides with the ramification divisor of $\pi_i$. 
We have in this way a collection $\{B_i\}$ of $\mathbf{Z}/a_i\mathbf{Z}$-stable 
Cartier divisors $B_i$ on $U_i$ such that $(p^n)_{i,j}^*(B_i) = (p^n)_{j,i}^*(B_j)$. 
Then $O(B) := \{O(B_i)\}$ becomes an orbifold line bundle and, by the definition of $B$, 
we have $O(B) \cong \mathcal{L}^{\otimes m}$ for some $m \geq 0$. 
Since $- K^{orb}_{\mathbf{P}(X)}$ and $O(B)$ are ample and nef respectively, the $\mathbf{Q}$-Cartier 
divisor $-K_{\mathbf{P}(X)}$ is ample.   
Since $\bar{\mathcal L}^{\otimes n} - K_{\mathbf{P}(X)}$ is ample and 
$(\mathbf{P}(X), 0)$ has log terminal singularities,  Lemma is a direct consequence of the Kawamata-Viehweg 
vanishing ([KMM], Theorem 1-2-5).
\vspace{0.2cm}

{\bf 6. Equivalence up to constant} 

Let $(X, \omega)$ be a pair of a normal affine variety $X$ of dimension $2d$ with a good $\mathbf{C}^*$-action and an algebraic symplectic 2-form $\omega$ on 
$X_{reg}$ with {\em positive} 
weight $l$. An algebraic symplectic 2-form $\omega'$ on $X_{reg}$ is said to 
be equivalent to $\omega$ up to constant when $\omega' = \lambda \cdot \omega$ 
with some $\lambda \in \mathbf{C}^*$. 
\vspace{0.2cm}

Let us consider the hypersurfaces 
$$ X_n:= \{(a,b,x,y,z) \in \mathbf{C}^5; a^2x + 2aby + b^2z + (xz-y^2)^n = 0\}, $$ 
where $n \geq 2$. These are central fibres of Slodowy slices to nilpotent orbits of 
$sp(2n)$ with Jordan type $[2n-2,1^2]$ ([LNS]); hence they admit natural symplectic 
2-forms $\omega_n$ of weight $2$. One can also construct symplectic 2-forms $\omega'_n$ 
of weight 2 on $X_n$ by using representations of $sl_2$ (cf. Introduction, Example 2, (b)). Moreover, 
$X_3$ coincides with the central fibre of the Slodowy slice to the subsubregular 
nilpotent orbit of the Lie algebra of type $G_2$ ([LNS], Section 10). Thus $X_3$ admits a symplectic 2-form 
$\sigma_3$ induced from the Kostant-Kirillov form on $\mathbf{g}_2$.  
By Theorem (3.1) we already know that they are equivalent up to $\mathbf{C}^*$-equivariant 
automorphism. But we can say more: 
\vspace{0.2cm}

{\bf Proposition (6.1)}. {\em Each hypersurface $X_n$ admits a unique holomorphic symplectic 2-form 
of weight $2$ up to constant.}
\vspace{0.2cm}

{\em Proof}.   
We put $X := X_n$. In this case, as explained below, 
$\mathrm{Codim}_{\mathbf{P}(X)}(\mathbf{P}(X) - \mathbf{P}(X)^0) = 2$ 
and $\mathbf{P}(X)^0 = \mathbf{P}(X)_{reg}$. 
As in (4.3), $\omega_n$ defines a contact form 
$\theta \in H^0(\mathbf{P}(X)_{reg}, \Omega^1_{\mathbf{P}(X)_{reg}}\otimes L^{\otimes 2})$. 
It is enough to check that $\theta$ is a unique contact structure with contact line bundle 
$L^{\otimes 2}$. 

First note that $\mathbf{P}(X)$ is not quasi-smooth; 
$X$ has a Du Val singularity of type $D_{n+1}$  
along $\{a = b = xz - y^2 = 0\}$. When $n = 2$, we understand that $D_3 = A_3$. 
The singular locus of $\mathbf{P}(X)$ is 
the disjoint union of two smooth rational curves 
$$ \{a = b = xz -y^2 = 0\} \cup \{x = y = z = 0\} $$ 
in $\mathbf{P}(2n-1,2n-1,2,2,2)$. Along the first component $\mathbf{P}(X)$ has 
a $D_{2n}$ surface singularity and along the second component it has 
quotient singularity of type $\frac{1}{2n-1}(1,1)$.    
Take points $p_1$ and $p_2$ respectively from the first and second components and 
consider the complex analytic germs $(\mathbf{P}(X), p_i)$. Then 
$$(\mathbf{P}(X), p_1) \cong (\mathbf{C}^1, 0) \times D_{2n}, $$  
$$(\mathbf{P}(X), p_2) \cong (\mathbf{C}^1, 0) \times \frac{1}{2n-1}(1,1). $$ 
Let $Cl(\mathbf{P}(X))$ (resp. $Cl(\mathbf{P}(X), p_i)$) be the divisor class 
group of $\mathbf{P}(X)$ (resp. $(\mathbf{P}(X), p_i)$). 
One has an exact sequence 
$$ 0 \to \mathrm{Pic}(\mathbf{P}(X)) \to  
Cl(\mathbf{P}(X)) \to \oplus_{1 \le i \le 2}Cl(\mathbf{P}(X), p_i).$$ 
By the same argument as in [Do, 3.2.5, 3.2.6], we see that 
$\mathrm{Pic}(\mathbf{P}(X)) = \mathbf{Z}\cdot [O_{\mathbf{P}(X)}(4n-2)]$. 
Since $Cl(\mathbf{P}(X), p_i)$ are finite abelian groups, 
we see that $Cl(\mathbf{P}(X))$ is a finitely generated Abelian group; in 
particular it is discrete. Let $\phi$ be an automorphism of $\mathbf{P}(X)$ 
contained in the neutral component $\mathrm{Aut}^0(\mathbf{P}(X))$ of 
the automorphism group $\mathrm{Aut}(\mathbf{P}(X))$. 
Then $\phi_*([O_{\mathbf{P}(X)}(i)]) = [O_{\mathbf{P}(X)}(i)]$ for 
all $i$.  
Note that there is an exact sequence 
$$ 0 \to O_{\mathbf{P}((2n-1)^2,2^3)}(i-4n) \to O_{\mathbf{P}((2n-1)^2,2^3)}(i) 
\to O_{\mathbf{P}(X)}(i) \to 0. $$ 
Applying these exact sequences, we have 
$$H^0(\mathbf{P}((2n-1)^2,2^3), O_{\mathbf{P}((2n-1)^2,2^3)}(i)) \cong 
H^0(\mathbf{P}(X), O_{\mathbf{P}(X)}(i))$$ 
for $i = 2, 2n-1$.  
Note that  $$H^0(\mathbf{P}((2n-1)^2,2^3), O_{\mathbf{P}((2n-1)^2,2^3)}(2)) 
= \mathbf{C}x \oplus \mathbf{C}y \oplus \mathbf{C}z,$$
and 
$$H^0(\mathbf{P}((2n-1)^2,2^3), O_{\mathbf{P}((2n-1)^2,2^3)}(2n-1)) 
= \mathbf{C}a \oplus \mathbf{C}b.$$
The automorphism $\phi$ induces linear automorphisms of 
$H^0(\mathbf{P}(X), O_{\mathbf{P}(X)}(i))$ $(i = 2, 2n-1)$  
and hence those of $\mathbf{C}x \oplus \mathbf{C}y \oplus \mathbf{C}z$ 
and $\mathbf{C}a \oplus \mathbf{C}b$. Such linear automorphisms induce an 
automorphism of $\mathbf{P}(2n-1, 2n-1, 2,2,2)$.  
Thus $\phi$ extends to  an automorphism of the ambient 
space $\mathbf{P}(2n-1, 2n-1, 2,2,2)$.   
  
We shall use Corollary (4.2.2) to prove the uniqueness of $\theta$. 
Let $j: \mathbf{P}(X)_{reg} \to \mathbf{P}(X)$ be the inclusion map. 
As we noted in (4.3), the contact structure $\theta$ induces an exact sequence  
$$ 0 \to j_*O(D) \to \Theta_{\mathbf{P}(X)} \to j_*(L^{\otimes 2}) \to 0.$$
Since $j_*(L^{\otimes 2}) = O_{\mathbf{P}(X)}(2)$, we see that 
$h^0(\mathbf{P}(X), j_*(L^{\otimes 2})) = 3$. 
 
On the other hand, $h^0(\mathbf{P}(X), \Theta_{\mathbf{P}(X)}) = 3$. 
A geometric explanation of this fact is the following. As we have seen above, all infinitesimal automorphisms 
of $\mathbf{P}(X)$ come from those of the ambient space $\mathbf{P}(2n-1,2n-1,2,2,2)$. 
The set of linear transformations of $(x,y,z)$ preserving the quadratic form $xz - y^2$ becomes  
a 3-dimensional algebraic subgroup of $GL(3, \mathbf{C})$. 
Fix such a linear transformation $\varphi$. Then there is a unique linear transformation 
of $(a,b)$ (up to sign) which sends the cubic form $a^2\varphi(x) + 2ab \varphi(y) + b^2 \varphi(z)$ 
to $a^2x + 2ab y + b^2 z$.  Since the exact sequence attached to the contact structure always 
splits (as $\mathbf{C}$-modules), we conclude that $h^0(j_*O(D)) = 0$. Q.E.D.   
\vspace{0.2cm}

Let $O \subset \mathbf{g}$ be a nilpotent adjoint orbit of a complex simple Lie algebra. 
Let $\tilde{O}$ be the normalisation of the closure $\bar{O}$. Since $O$ admits a Kostant-Kirillov 
2-form, $\tilde{O}$ has a holomorphic symplectic structure of weight $1$. 

{\bf Proposition (6.2)} {\em Assume that $\tilde{O}$ is a Richardson orbit with a Springer map $\pi: T^*(G/P) \to 
\tilde{O}$ for some parabolic subgroup $P$ of $G$. Then $\tilde{O}$ has a unique symplectic 
structure of weight $1$ up to constant.}

{\em Proof}.  Let $\mathbf{P} := \mathbf{P}(T^*(G/P))$ be the projectivised tangent bundle of $G/P$. 
Then $\pi$ induces a generically finite proper map $\bar{\pi}: \mathbf{P} \to \mathbf{P}(\bar{O})$ 
and the contact 1-form $\theta \in H^0(\mathbf{P}(O), \Omega^1_{\mathbf{P}(O)}\otimes O_{\mathbf{P}(O)}(1))$ 
is pulled back (and is extended) to a contact 1-form $$\bar{\pi}^*\theta \in 
H^0(\mathbf{P}, \Omega^1_{\mathbf{P}} \otimes O_{\mathbf{P}}(1)).$$  
We prove that this is a unique contact structure on $\mathbf{P}$ with contact line bundle 
$O_{\mathbf{P}}(1)$. Let 
$$ 0 \to O(D) \to \Theta_{\mathbf{P}} \stackrel{\bar{\pi}^*\theta}\to O_{\mathbf{P}}(1) \to 0 $$ 
be the corresponding exact sequence. Let $p: \mathbf{P} \to G/P$ be the projection map 
of the projective space bundle.  Since $p_*O_{\mathbf{P}}(1) = \Theta_{G/P}$, we have 
$$h^0({\mathbf P}, O_{\mathbf P}(1)) = h^0(G/P, \Theta_{G/P}).$$   
On the other hand, by the exact sequences  
$$ 0 \to O_{\mathbf P} \to p^*\Omega^1_{G/P}\otimes O_{\mathbf P}(1) \to \Theta_{{\mathbf P}/(G/P)} 
\to 0, $$ one has an exact sequence 
$$ 0 \to H^0(O_{\mathbf P}) \to H^0(\underline{\mathrm Hom}(\Theta_{G/P}, \Theta_{G/P}) 
\to H^0(\Theta_{{\mathbf P}/(G/P)}) \to H^1(O_{\mathbf{P}}) .$$ 
Since $\Theta_{G/P}$ is a simple vector bundle ([A-B]), we have   
$H^0(\underline{\mathrm Hom}(\Theta_{G/P}, \Theta_{G/P}) \cong \mathbf{C}$.  
As $H^1(O_{\mathbf{P}}) = 0$,  we  see that $H^0(\Theta_{{\mathbf P}/(G/P)}) = 0$. 
By the exact sequence 
$$ 0 \to H^0(\Theta_{{\mathbf P}/(G/P)}) \to H^0(\Theta_{\mathbf P}) \to 
H^0(p^*\Theta_{G/P}), $$ it is clear that 
$h^0(\Theta_{\mathbf P}) = h^0(G/P, \Theta_{G/P})$.  This implies that $H^0(\mathbf{P}, O(D)) = 0$. 
Q.E.D. 
\vspace{0.2cm}

{\bf Remark}.  Let  $O$ be a nilpotent orbit (where $O$ is not necessarily  a Richardson orbit).  
Consider the contact structure on $\mathbf{P}(O)$: 
$$ 0 \to O(D) \to \Theta_{{\mathbf P}(O)} \stackrel{\theta}\to O_{{\mathbf P}(O)}(1) \to 0.$$  
Since $O$ is a homogeneous space acted by $G$, there 
is a natural map $\mathbf{g} \to H^0(\Theta_{{\mathbf P}(O)})$. Then the composition map 
$$ \theta\vert_{\mathbf g}: \mathbf{g} \to H^0(O_{{\mathbf P}(O)}(1)) $$ is injective. 
The following is a proof. Let $\omega$ be the Kostant-Kirillov 2-form on $O$. 
As in (4.3), let $\zeta$ be the vector field on $O$ which generates the $\mathbf{C}^*$-action. 
Let $\pi: O \to \mathbf{P}(O)$ be the projection map. By definition, $\pi^*\theta = 
\omega(\zeta, \cdot)$. For $x \in O$, we denote by $\bar{x} \in \mathbf{P}(O)$ the corresponding point.
Let us consider $T_xO$ as a linear subspace of ${\mathbf g}$.  
Then $\zeta_x = x$ by the definition.  For $v \in \mathbf{g}$,  we have $[x, v] \in T_xO$; hence 
$$(\theta\vert_{\mathbf g}(v))_{\bar{x}} = \omega_x(x, [v,x]).$$ 
One can write $x = [a_x, x]$ with some $a_x \in \mathbf{g}$. 
Let $\kappa$ be the Killing form on $\mathbf{g}$. By the definition of the Kostant-Kirillov 
2-form we have 
$$  \omega_x(x, [v,x]) = \kappa (x, [a_x, v]) = \kappa ([x, a_x], v) = -\kappa (x, v).$$ 
If $v \in \mathrm{ker}(\theta_{\mathbf g})$, then $\kappa (x, v) = 0$ for all $x \in O$. 
Note that, $x$ is contained in the cone $\bar{O} \subset \mathbf{g}$. 
Since $T_0{\bar{O}}$ is invariant under the adjoint $G$-action and the adjoint representation is irreducible, 
$T_0{\bar{O}} = \mathbf{g}$. This means that, if $x$ runs inside $O$, they span $\mathbf{g}$ as 
a $\mathbf{C}$-vector space. Since $\kappa$ is non-degenerate, we conclude that $v = 0$.     
Now we have 
\vspace{0.2cm}

{\bf Problem}. {\em When does ${\mathbf g}$ coincide with $H^0(\mathbf{P}(O), \Theta_{\mathbf{P}(O)})$ ?} 
\vspace{0.2cm}

When $O_{min}$ is the minimal nilpotent orbit of {\bf g}, $\mathbf{P}(O_{min})$ is a flag variety $G/P$ 
with a parabolic subgroup $P$.  Let  
$M := G/P$ be a flag variety  where $G$ is a connected simple complex 
Lie group acting effectively on $M$. Then, by Onishchik (cf. [G-O, Theorem 4.10]),   
the neutral component $\mathrm{Aut}^0(G/P)$ is isomorphic to $G$ except in the following three 
cases. 
\vspace{0.2cm}

(i) $G = PSp(2n)$ and $P$ is the stabilizer subgroup of an isotropic flag of type 
$(1, 2n-2, 1)$ in the vector space $\mathbf{C}^{2n}$ acted by $G$. 
\vspace{0.2cm}

(ii) $G = G_2 \subset SO(7)$ and  $M$ is a quadric 5-fold in $\mathbf{P}^6$. 
\vspace{0.2cm}

(iii) $G = SO(2n+1)$ and $P$ is the stabilizer subgroup of an isotropic 
flag of type $(n,1,n)$ in $\mathbf{C}^{2n+1}$.  
\vspace{0.2cm}

In (ii) and (iii), $M = G/P$ is not realized as the projectivised cone 
$\mathbf{P}(O_{min})$ of the minimal nilpotent orbit $O_{min}$.   
But in the case (i), $G/P = \mathbf{P}(O_{min})$ with $O_{min} \subset sp(2n)$. 
Thus we have proved the following.  
\vspace{0.2cm}

{\bf Proposition (6.3)} {\em Assume that $O_{min}$ is the minimal nilpotent orbit of $\mathbf{g}$. 
Then $\tilde{O}_{min}$ has a unique symplectic structure of weight $1$ up to constant except 
when $\mathbf{g} = sp(2n)$.} 
\vspace{0.2cm}

Note that the exceptional case corresponds to the quotient singularity 
$\mathbf{C}^{2n}/\mathbf{Z}_2$ by the action $(z_1, ..., z_{2n}) \to 
(-z_1, ..., -z_{2n})$. 
\vspace{0.2cm}

{\bf 7. Problems} 

Let $(X, \omega)$ be a pair of a normal affine variety $X$ of dimension $2d$ with a good $\mathbf{C}^*$-action and an algebraic symplectic 2-form $\omega$with 
{\em positive} weight $l$. 
Let us call $(X, \omega)$ {\em irreducible} of weight $l$ when 
$\omega$ is a unique symplectic structure of weight $l$ up to 
constant. 

{\bf Problem (7.1)}: {\em Does $(X, \omega)$ have symplectic singularities, or equivalently, 
canonical singularities ?}
\vspace{0.15cm}

{\bf Problem (7.2)}: {\em Is the fundamental group $\pi_1(X_{reg})$ of the regular part 
of $X$ finite ?} 
\vspace{0.2cm}

When $G:= \pi_1(X_{reg})$ is finite, one can take a finite $G$-Galois covering 
$\pi: Y \to X$ in such a way that the induced map $\pi^{-1}(X_{reg}) \to X_{reg}$ 
is the universal covering of $X_{reg}$. Let $m$ be the order of $G$. 
Let $\mathbf{C}^* \times X \to X$ $(t, x) \to \phi_t(x)$ be the given $\mathbf{C}^*$-action 
on $X$. We consider the $\mathbf{C}^*$-action on $X$ defined as its $m$-th power: 
$$ \mathbf{C}^* \times X \to X\;  (t, x) \to \phi_{t^m}(x).$$ 
Then $Y$ has a $\mathbf{C}^*$-action so that $\pi$ is $\mathbf{C}^*$-equivariant. 

Recall here the Bogomolov splitting theorem for a compact K\"{a}hler manifold $X$ with $c_1 = 0$. 
It states, in particular, that if $X$ is a holomorphic symplectic manifold with a finite fundamental group, 
then its universal cover $\tilde{X}$ splits into the product of irreducible symplectic manifolds 
$X_i$ ($i = 1, ..., r$) such that $h^0(X_i, \Omega^2_{X_i}) = 1$.

The following is an analogue of the splitting theorem in affine symplectic varieties 
with good $\mathbf{C}^*$-actions.  
\vspace{0.2cm}

{\bf Problem (7.3)}: {\em Is there a $\mathbf{C}^*$-equivariant isomorphism of 
symplectic varieties} 
$$(Y, \pi^*\omega) \cong \prod_{1 \le i \le k}(Y_i, \omega_i)$$ 
{\em where each $(Y_i, \omega_i)$ is irreducible of weight $m\cdot l$ ?} 
\vspace{0.2cm}

For example, as an $(X, \omega)$, take the quotient singularity $\mathbf{C}^{2n}/\mathbf{Z}_2$ 
defined at the end of {\bf 6} and the symplectic form induced from 
$\tilde{\omega} := dz_1 \wedge dz_2 + ... + dz_{2n-1}\wedge dz_{2n}$. Then $(X, \omega)$ is not irreducible; but 
$$(\mathbf{C}^{2n}, \tilde{\omega}) \cong \prod_{1 \le i \le n} (\mathbf{C}^2, dz_{2i-1}\wedge dz_{2i}).$$

{\bf Acknowledgement}. I would like to thank M. Lehn for discussions on Lemma (1.2).  

\vspace{0.2cm}

\begin{center}
{\bf References} 
\end{center}

[A-B] Azad, H., Biswas, I.: A note on the tangent bundle of G/P, Proc. Indian Acad. Sci. {\bf 120}, No 1 
(2010), 69-71

[Ba] Baily, W.L.: On the embedding of V-manifolds in projective spaces, Amer. J. Math. 
{\bf 79} (1957), 403-430 

[Be] Beauville, A.: Symplectic singularities, Invent. Math. {\bf 139} (2000), 541-549

[Do] Dolgachev, I.: Weighted projective varieties, Group actions and vector fields. Lecture 
Notes in Mathematics {\bf 956} (1982) 

[C-M] Collingwood, D.H., McGovern, W.M.: Nilpotent orbits in semisimple Lie algebras, 
Van Nostrand Reinhold Math. Ser. Van Reinhold, New York, 1993 

[E-G] Etingof, P., Ginzburg, V.: Noncommutative del Pezzo surfaces and 
Calabi-Yau algebras, J. Eur. Math. Soc. {\bf 12} (2010) 1371-1416

[Fl] Flenner, Extendability of differential forms on nonisolated singularities, 
Invent. Math. {\bf 94} (1988), 317-326 

[G-O] Gorbatsevich, V., Onishchik, A.:  Lie transformation groups, 95-229: in Lie groups 
and Lie algebras I, Encyclopaedia of Mathematical Sciences {\bf 20}, Springer-Verlag  

[H-M] Hanany, A., Mekareeya, N.: Tri-vertices and SU(2)'s. arxiv:1012.2119v1. 

[Hi] Hironaka, H.: Resolution of singularities of an algebraic variety over 
a field of characteristic zero, Ann. Math. {\bf 79} (1964), 109-326

[Ka] Kaledin, D.: Symplectic singularities from the Poisson point of view, 
J. Reine Angew. Math. {\bf 600} (2006), 135-156   

[Kaw] Kawamata, Y.: The cone of curves of algebraic varieties, 
Ann. of Math. {\bf 119} (1984) no 3. 603-633

[KMM] Kawamata, Y., Matsuda, K., Matsuki, K.: 
Introduction to the minimal model program, Adv. Pure Math. {\bf 10}, 
283-360, Kinokuniya-North-Holland, 1987

[Ko] Koll\'{a}r, J.: Singularities of pairs, Proc. Symp. Pure Math. {\bf 62}-1, (1997) 
221-287

[LeB] LeBrun, C.: Fano manifolds, contact structures, and quartaionic geometry, 
International Journal of Mathematics {\bf 6}, No 3 (1995) 419-437 

[Li] Lichnerowicz, A.: Varietes de Jacobi et spaces homogenes de contact complexes, 
J. Math. Pure et Appl. {\bf 67} (1988) 131-173

[LNS] Lehn, M., Namikawa, Y., Sorger, C.: Slodowy slices and universal Poisson deformations, 
to appear in Compositio Math.  

[L-N-S-vS] Lehn, M., Namikawa, Y., Sorger, C., van Straten,D.: Symplectic 
hypersurfaces, in preparation  

[Mo] Moser, J.: On the volume elements on a manifold, Trans. Amer. Math. Soc. {\bf 120} 
(1965), 286-294  

[Mu] Mumford, D.: Towards an enumerative geometry of moduli spaces of curves, in Arithmetic 
and Geometry, edited by M. Artin, J. Tate, Birkhauser-Boston, 1983, 271-326

[Na 1] Namikawa, Y.: Poisson deformations of affine symplectic varieties, Duke Math. J. {\bf 156} 
(2011) 51-85 

[Na 2] Namikawa, Y.: Extension of 2-forms and symplectic varieties, J. Reine Angew. Math. {\bf 539} 
(2001), 123-147

[Na 3] Namikawa, Y.: Deformation theory of singular symplectic n-folds, Math. Ann. {\bf 319} (2001), 
597-623

[Naru] Naruki, I.: Some remarks on isolated singularity and their application to algebraic manifolds.  
Publ. Res. Inst. Math. Sci.  {\bf 13}  (1977/78), no. 1, 17-46

[Sl] Slodowy, P.: Simple singularities and simple algebraic groups, Lecture Notes in 
Mathematics, {\bf 815} (1980), Springer Verlag 

\vspace{0.2cm}

\begin{center}
Department of Mathematics, Faculty of Science, Kyoto University, Japan 

namikawa@math.kyoto-u.ac.jp
\end{center}

\end{document}